\documentstyle[amssymb,12pt]{article}
\title{$\mathbf P$-NDOP and $\mathbf P$-decompositions of 
$\aleph_\epsilon$-saturated models of superstable
theories}
\author{Saharon Shelah\thanks{Partially supported by the 
US-Israel Binational Science Foundation 
and by NSF grants DMS-0600940 and DMS-1101597. Publication 933.}\\
The Hebrew University of Jerusalem\\
Einstein Institute of Mathematics\\
Edmond J.\ Safra Campus, Givat Ram\\
Jerusalem, 91904, Israel \\ \\
Department of Mathematics\\
Hill Center-Busch Campus\\
Rutgers, the State University of New Jersey\\
110 Frelinghuysen Road\\Piscataway, NJ 08854-9019 USA
\and Michael C.\
Laskowski\thanks{Partially supported
by NSF grant DMS-0901336.}\\
Department of Mathematics\\University of Maryland\\
College Park, MD 20742 USA
}

\newbox\smilebox
\newbox\anchorbox
\newbox\noanchorbox
\newbox\tempbox

\setbox\smilebox=\hbox{$\smile$}

\def\anchor{\hbox{\vtop{
           \hbox to \wd\smilebox{\hfil\vrule width.4pt height7pt depth1pt\hfil}
           \vskip  -11.5truept
           \hbox to \wd\smilebox{\hfil$\smile$\hfil}}}}
\setbox\anchorbox=\anchor
\def\noanchor{\hbox{\vtop{
           \hbox to \wd\anchorbox{\hfil\anchor\hfil}
           \vskip -14truept
           \hbox to \wd\anchorbox{\hfil/\hfil}}}}
\setbox\noanchorbox=\noanchor

\def\fg#1#2#3{\setbox\tempbox=\hbox{$\scriptstyle{#2}$}
\ifnum\wd\anchorbox>\wd\tempbox\dimen255=\wd\anchorbox
\else\dimen255=\wd\tempbox\fi
{#1\,\vtop{\hbox to \dimen255{\hfil\anchor\hfil}
           \vskip -6truept
           \hbox to \dimen255{\hfil$\scriptstyle{#2}$\hfil}}
           \,#3}}

\def\nfg#1#2#3{\setbox\tempbox=\hbox{$\scriptstyle{#2}$}
\ifnum\wd\noanchorbox>\wd\tempbox\dimen255=\wd\noanchorbox
\else\dimen255=\wd\tempbox\fi
{#1\,\vtop{\hbox to \dimen255{\hfil\noanchor\hfil}
           \vskip -6truept
           \hbox to \dimen255{\hfil$\scriptstyle{#2}$\hfil}}
           \,#3}}

\setbox1=\hbox{$\bot$}

\def\north#1#2{#1\,
\hbox{$\bot$\llap {\hbox to\wd1 {\hfil $/$\hfil}}}
\,#2}

\def\nao#1#2#3{#1\  \hbox{\vtop{ 
\baselineskip=4pt
\hbox{$\bot$\llap {\hbox to\wd1 {\hfil $/$\hfil}}
\hskip .05em \llap{\hbox{$^{\scriptscriptstyle{a}}$}}}\hbox{$\scriptstyle
{#2}$}}}\, #3}

\def\bp{\par{\bf Proof.}$\ \ $}

\def\includeE#1{{\lhook\kern-3.5pt\joinrel\smash{
    \mathop{\longrightarrow}\limits^{#1}}}}

\def\efor/{Example~\ref{E4}}

\def\BL/{Baldwin--Lachlan}
\def\Bu/{Buechler}
\def\Hr/{Hrushovski}
\def\lm/{locally modular}
\def\wm/{weakly minimal}
\def\nm/{non--modular}
\def\ss/{superstable}
\def\ud/{unidimensional}
\def\sm/{strongly minimal}

\def\abar{\overline{a}}

\def\bbar{\overline{b}}

\def\hbar{\overline{h}}

\def\lbar{\overline{l}}

\def\Mbar{\overline{M}}

\def\acl{{\rm acl}}

\def\dom{{\rm dom}}

\def\tp{{\rm tp}}
\def\stp{{\rm stp}}

\def\tr/{trivial}
\def\nt/{non--trivial}
\def\st/{strong type}

\def\TV/{Tarski--Vaught}

\def\sc/{sound construction}
\def\ac/{atomic construction}

\def\fal/{functional}
\def\upl/{unique parallel lines}
\def\chp/{categorical in a higher power}

\def\conc{{\char'136}}
\def\abar{\bar{a}}
\def\bbar{\bar{b}}

\def\phi{\varphi}

\def\C{{\frak  C}}

\def\tp{{\rm tp}}
\def\stp{{\rm stp}}

\def\contains{\supseteq}
\def\dom{{\rm dom}}
\def\acl{{\rm acl}}
\def\dcl{{\rm dcl}}
\def\bp{{\bf Proof.}\quad}
\def\endproof{\medskip}
\def\<{\langle}
\def\>{\rangle}
\def\d{{\mathfrak d}}
\def\P{{\mathbf P}}
\def\Pr{{\mathbf P^r}}
\def\Pperp{\P^{\perp}}
\def\Pactive{\P^{{\mathbf active}}}
\def\Pdull{\P^{{\mathbf dull}}}
\def\Treg{{\bf T}^{\rm reg}}
\def\lbar{\overline{\lambda}}
\def\lbarP{(\lbar,\P)}
\def\eps{\epsilon}
\def\ae{\aleph_\eps}
\def\aeP{(\ae,\P)}
\def\PP{{\cal P}_\P}
\def\dpP{{\rm dp}_{\P}}
\newtheorem{Theorem}{Theorem}[section]
\newtheorem{Proposition}[Theorem]{Proposition}
\newtheorem{Definition}[Theorem]{Definition}
\newtheorem{Convention}[Theorem]{Convention}

\newtheorem{Lemma}[Theorem]{Lemma}
\newtheorem{Corollary}[Theorem]{Corollary}
\newtheorem{Fact}[Theorem]{Fact}
\begin{document}
\maketitle
\begin{abstract}
Given a complete, superstable theory, we distinguish a 
class $\P$  of regular types, typically closed under automorphisms
of $\C$ and non-orthogonality.  We define the notion of $\P$-NDOP,
which is a weakening of NDOP.
For superstable theories with $\P$-NDOP, we
prove the existence of $\P$-decompositions
and derive an analog of \cite{Sh401}.
In this context, we also find a sufficient condition on
$\P$-decompositions that imply non-isomorphic models.
For this, we investigate natural structures on the types in $\P\cap S(M)$
modulo non-orthogonality.   
\end{abstract}

\section{Introduction}

Results by the first author, most notably Chapter~X of \cite{Shc}
and the first half of \cite{Sh401} demonstrate that $\aleph_\epsilon$-saturated models
of superstable theories with NDOP admit very desirable decompositions.
In this paper, we generalize these results in three ways.
First, we always assume that the theory
$T$ is superstable, but we only have NDOP for a class $\P$ of regular types.
Second, we show that the tree structure of a decomposition of
an $\aleph_\eps$-saturated model $M$ can be read off
from the non-orthogonality classes of regular types in $S(M)$.  Third,
we show that these results for $\aleph_\eps$-saturated models
give information about weak decompositions of arbitrary models of 
such theories.  

In more detail, throughout the paper we assume we have a fixed,
complete, {\bf superstable} theory and we work within a monster model
$\C$.  We fix a set $\P$ of stationary, regular types over small subsets of
$\C$ that is closed under automorphisms of $\C$ and the equivalence relation of
nonorthogonality, and additionally assume that our theory satisfies $\P$-NDOP.
Typically, we fix a model $M$ that is at least $\aleph_\eps$-saturated
(i.e., $M$ contains a realization of every strong type over every finite subset of $M$)
and study $\P$-decompositions inside $M$ of many varieties.  Of primary interest
are prime, $\aeP$-decompositions $\d$  of $M$ over $B\choose A$ 
(see Definition~\ref{ofMoverBA})
where $A\subseteq B$ are $\eps$-finite and every regular type $p$ non-orthogonal to
$\stp(B/A)$ is in $\P$.  We associate a subset $\PP(\d,M)$ of $S(M)\cap\P$
(see Definition~\ref{PP}) to such a pair.  
The main theorem of the paper, Theorem~\ref{big}, asserts that this set of regular types
depends only on $B\choose A$.
In particular, it is independent of the decomposition $\d$, and successive results show that these sets
have a tree structure under inclusion.  

In the final section of the paper, we show how this result, which holds only for $\ae$-saturated
models, gives positive information for much weaker decompositions of models $M_0$ without any
saturation assumption.  

\section{Preliminaries}

As mentioned above, we always work in a class of models of a complete, superstable, first-order theory
$T$.  We fix a monster model $\C$, and all models and sets we discuss will be small subsets of $\C$.
We assume that $T$ eliminates quantifiers, so any model $M$ will be an elementary submodel of $\C$,
and we additionally assume that `$T=T^{{\rm eq}}$', so that every type over an algebraically closed
set is stationary.

\begin{Definition} {\em   A set $A$ is {\em $\eps$-finite\/}
if $\acl(A)=\acl(a)$ for some $a\in\C^{\rm eq}$. 
}
\end{Definition}

Recall that as we are working in $\C^{\rm eq}$, it would be equivalent
to say that $\acl(A)=\acl(\abar)$ for some finite tuple.
It is easily seen that the union of two $\eps$-finite sets is $\eps$-finite.
Furthermore, since $T$ is superstable, any subset $B\subseteq A$ of
an $\eps$-finite set is $\eps$-finite.  [Why?  If $B\subseteq A$ with 
$\acl(A)=\acl(a)$,  choose a finite $\bbar$ from $B$ such that
$\fg B {\bbar} a$.  Then $\acl(B)=\acl(\bbar)$.]
Thus, the set of $\eps$-finite subsets of $\C$ form an ideal.

\begin{Convention} {\em  $\aleph_\eps$ is a cardinal strictly between
$\aleph_0$ and $\aleph_1$.
}
\end{Convention}

Thus, if we write `$M$ is $\lambda$-saturated for some $\lambda\ge\aleph_\eps$'
we mean that either $M$ is $\aleph_\eps$-saturated (i.e., realizes all  types
over $\eps$-finite subsets) or $M$ is $\lambda$-saturated for some $\lambda\ge\aleph_1$.
Recall that by e.g., IV~2.2(7) of \cite{Shc}, that for $\lambda\ge\aleph_1$, $M$
is $\lambda$-saturated if and only if  for every subset $A\subseteq M$ of
size less than $\lambda$, $M$ realizes every  type over $\acl(A)$.

We record several facts from \cite{Shc} that will be used throughout this paper.
The first is the Second Characterization Theorem, IV~4.18, the second is
X~Claim~1.6(5), the third is V~1.12, and (4) follows easily from (2) and (3).

\begin{Fact}  \label{refer}
Suppose $T$ is superstable and  $\lambda\ge\aleph_\eps$. 
\begin{enumerate}
\item A  model $M$ is $\lambda$-prime over a set $A$
if and only if (1) $M\contains A$ and is $\lambda$-saturated; (2) $M$ is $\lambda$-atomic over
$A$; and (3) every $A$-indiscernible sequence $I\subseteq M$ has length at most $\lambda$.
(When $\lambda=\aleph_\eps$, the $\lambda$ occurring in (3) should be replaced by $\aleph_0$.)

\item  If $M$ is $\lambda$-saturated, $A\contains M$, and $N$ is $\lambda$-prime over $M\cup A$,
then $N$ is dominated by $A$ over $M$.

\item  If $M\preceq N$ are both $\lambda$-saturated, $p\in S(M)$ is regular,
and there is some $c\in N\setminus M$ such that $\tp(c/M)\not\perp p$,
then $p$ is realized in $N$.

\item  If $M_0\preceq M_1\preceq M_2$ are all $\lambda$-saturated
and there is $e\in M_2\setminus M_1$ such that $\tp(e/M_1)$ is regular
and non-orthogonal to $M_0$, then there is $e^*\in M_2\setminus M_1$
such that $e$ and $e^*$ are domination equivalent over $M_1$,
with $\fg {e^*} {M_0}  {M_1}$.
\end{enumerate}
\end{Fact}

\section{P-NDOP}

Our story begins by localizing the notion of DOP around a single
parallelism class of stationary, regular types.

\begin{Definition} \label{1.1}
{\em   An {\em independent triple\/}
of models $(M_0,M_1,M_2)$ satisfy $M_0=M_1\cap M_2$ and
$\{M_1,M_2\}$ are independent over $M_0$.
For $\lambda\ge\aleph_{\epsilon}$,
a {\em $\lambda$-quadruple} is a sequence $(M_0,M_1,M_2,M_3)$
of $\lambda$-saturated models, where $(M_0,M_1,M_2)$ form
an independent triple, and $M_3$ is $\lambda$-prime over
$M_1\cup M_2$.   A {\em $\lambda$-DOP witness\/} for a
stationary, regular type $p$ is a $\lambda$-quadruple
$(M_0,M_1,M_2,M_3)$ for which $Cb(p)\subseteq M_3$,
but $p\perp M_1$ and $p\perp M_2$.  We say that {\em $p$ has
a DOP witness\/} if it has a $\lambda$-DOP witness for some
$\lambda\ge\aleph_\eps$. 
}
\end{Definition}

Visibly,  whether a specific $\lambda$-quadruple is a $\lambda$-DOP
witness for $p$ depends only on the parallelism class of $p$.
To understand the consequences of this notion, we recall that 
a set $A$ is {\em self-based} on an independent triple
$(M_0,M_1,M_2)$ of models if $\fg A {A\cap M_i} {M_i}$ holds for
each $i<3$.  The concept of self-basedness was defined explicitly
in \cite{LSh} and was used implicitly in the proof of X~2.2(iii$\rightarrow$iv)
of \cite{Shc}.   The fact that for  any independent triple $(M_0,M_1,M_2)$, any finite set $A$ can be extended to a  finite, self-based $B\subseteq AM_1M_2$ 
follows from Lemma~2.4 of \cite{LSh}.  
The Claim in the proof of Theorem~1.3 of \cite{LSh}
 establishes the following Fact.

\begin{Fact}  \label{selfbased} 
If $A$ is self-based on the independent triple
$(M_0,M_1,M_2)$, $p\in S(A)$ is stationary, $p\perp M_1$, and
$p\perp M_2$, then $p\vdash p|AM_1M_2$.
\end{Fact}

Using this Fact,
an easy examination of the proof of \cite{Shc}, X~2.2 yields:

\begin{Fact}  \label{X2.2}
Let $p$ be any  stationary, regular type with a DOP witness.  Then:
\begin{enumerate}
\item  For every $\lambda\ge\aleph_\eps$, $p$ has a $\lambda$-DOP witness;
\item  For every $\lambda$-DOP witness $(M_0,M_1,M_2,M_3)$
for $p$, there is an infinite,
indiscernible set $I\subseteq M_3$ over $M_1\cup M_2$ whose average
type $Av(I,M_3)$ is parallel to $p$; and
\item  For every $\lambda$-DOP witness $(M_0,M_1,M_2,M_3)$
for $p$, there is a subset $A\subseteq M_3$, $|A|<\lambda$ over
which $p$ is based and stationary and a Morley sequence $\<b_i:i<\lambda\>$
from $M_3$ in $p|AM_1M_2$.
\end{enumerate}
\end{Fact}

We isolate one Corollary from this that will be crucial for us later.

\begin{Corollary}  \label{isom}
For any $\lambda\ge\aleph_\eps$, if $(M_0,M_1,M_2,M_3)$
is a $\lambda$-DOP witness for a stationary, regular $p\in S(M_3)$,
then for any realization $c$ of $p$,  any $\lambda$-prime model
$M_3[c]$ over $M_3\cup\{c\}$ is isomorphic to $M_3$ over
$M_1\cup M_2$.  In particular, $M_3[c]$ is $\lambda$-prime over 
$M_1\cup M_2$.
\end{Corollary}

\bp  By the uniqueness of $\lambda$-prime models,
both statements will follow once we establish that $M_3\cup\{c\}$
is the universe of a $\lambda$-construction sequence over $M_1\cup M_2$.
To see this, first fix a $\lambda$-construction sequence $\<b_i:i<\delta\>$
of $M_3$ over $M_1\cup M_2$.  As notation, for each $i<\delta$,
let $B_i=M_1\cup M_2\cup\{b_j:j<i\}$ and fix a subset 
$X_i\subseteq B_i$, $|X_i|<\lambda$ such that $\stp(b_i/X_i)\vdash\stp(b_i/B_i)$.

Next, choose a subset $A\subseteq M_3$, $|A|<\lambda$
 over which $p$ is based and stationary.    By forming an increasing
$\omega$-chain, we can increase $A$ slightly (still maintaining $|A|<\lambda$)
so that $A$ is self-based on $(M_0,M_1,M_2)$ and $X_i\subseteq A$ whenever
$b_i\in A$.

Let $\<a_i:i<\gamma\>$ be the enumeration of $A$ given by the ordering
of the original construction.  Easily, $\<a_i:i<\gamma\>$ is
 $\lambda$-constructible
over $M_1\cup M_2$.  

Furthermore, it follows from Fact~\ref{selfbased} that for any Morley
sequence $I$ in $p|A$ with $|I|<\lambda$, we have $p|AI\vdash p|AIM_1M_2$.
Using this, we have a $\lambda$-construction sequence 
$\<a_i:i<\gamma\>\conc\<c_j:j<\lambda\>$ over $M_1\cup M_2$,
where $\<c_j:j<\lambda\>$ is any Morley sequence in $p|A$
from $M_3$ (the existence of such a sequence follows from Fact~\ref{X2.2}(4)).
It follows from the uniqueness of $\lambda$-prime models and the fact
that such models are $\lambda$-constructible that there is another
$\lambda$-construction sequence of $M_3$ over $M_1\cup M_2$
in which 
$\<a_i:i<\gamma\>\conc\<c_j:j<\lambda\>$ is an initial segment.
As notation, let $\<b_k:k<\nu\>$ be the tail of this sequence.
For each $k<\nu$, let 
$B_k^*=M_1\cup M_2\cup A\cup\{c_j:j<\lambda\}\cup \{b_\ell:\ell<k\}$
and  choose $Y_k\subseteq B_k^*$, $|Y_k|<\lambda$ such that
$\stp(b_k/Y_k)\vdash \stp(b_k/B_k^*)$.  Without loss, we may assume
$A\subset Y_k$ for each $k$.
To complete the proof, it suffices to prove that 
$$\<a_i:i<\gamma\>\conc\<c\>\conc\<c_i:i<\lambda\>\conc\<b_k:k<\nu\>$$
is a $\lambda$-construction sequence over $M_1\cup M_2$.  

We already know that $\<a_i:i<\gamma\>$ is a $\lambda$-construction
sequence over $M_1\cup M_2$.
Using the first sentence of the previous paragraph, combined with the fact
that $\{c\}\cup\{c_j:j<\lambda\}$ is independent over $A$, we inductively
obtain that $\<a_i:i<\gamma\>\conc\<c\>\conc\<c_j:j<\lambda\>$
is also a  $\lambda$-construction sequence over $M_1\cup M_2$.
Thus, it suffices to prove that $\stp(b_k/Y_k)\vdash\stp(b_k/B^*_kc)$
for each $k<\nu$.  For this, since both  $\tp(c/B^*_k)$
and $\tp(b_k/B^*_k)$ do not fork over $Y_k$, it suffices to show that
$\tp(c/Y_k)$ is almost orthogonal to $\stp(b_k/Y_k)$.  
To see this, choose $j<\lambda$ such that $\tp(c_j/A)$ does not fork over
$Y_k$.  Now, $\tp(c/Y_k)=\tp(c_j/Y_k)$ and $\tp(c_j/Y_k)$ is
almost orthogonal to $\stp(b_k/Y_k)$ since $\stp(b_k/Y_k)\vdash\stp(b_k/Y_kc_j)$,
so we finish.
\endproof

Next, we show additional closure properties of DOP witnesses.

\begin{Definition}{\em  A regular type $q$ {\em lies directly above $p$\/}
if there is a non-forking extension $p'\in S(M)$ of $p$ with $M$
$\aleph_\eps$-saturated, a realization $c$ of $p'$, and an
$\aleph_\eps$-prime model $M[c]$ over $M\cup\{c\}$ such that
$q\not\perp M[c]$, but $q\perp M$.
A regular type $q$ {\em lies above $p$\/} if there is a sequence
$p_0,\dots,p_n$ of types such that $p_0=p$, $p_n=q$, and
$p_{i+1}$ lies directly above $p_i$ for each $i<n$.  (We allow $n=0$,
so in particular, any regular type lies above itself.)

We say that $p$ {\em supports} $q$ if $q$ lies above $p$.
}
\end{Definition}

The nomenclature above is apt if one considers a branch of a decomposition tree.
Suppose $M_0\preceq\dots \preceq M_n$ is a sequence of $\aleph_\eps$-saturated
models such that for each $i<n$ there is $a_i\in M_{i+1}$ such that
$\tp(a_i/M_i)$ is regular (and orthogonal to $M_{i-1}$ when $i>0$) and
$M_{i+1}$ is $\aleph_\eps$-prime over $M_i\cup\{a_i\}$.  Then any
regular $q\not\perp M_n$ lies over any regular type $p$ non-orthogonal to
$\tp(a_0/M_0)$.  Similarly, any such $p$ supports any such $q$.

\begin{Proposition} \label{DOPclosure}
Fix a stationary, regular type $p$ with a  DOP witness.  Then:
\begin{enumerate}
\item  Every type parallel to $p$ has a DOP witness;
\item  Every automorphic image of $p$ has a DOP witness;
\item  Every stationary, regular $q$ non-orthogonal to $p$ has a DOP witness;
\item  Every stationary, regular  $q$
lying above $p$ has a DOP witness.
\end{enumerate}
\end{Proposition}

\bp  (1) and (2) are immediate.  
For (3), choose $\lambda\ge\aleph_\eps$ and a $\lambda$-quadruple
$(M_0,M_1,M_2,M_3)$ witnessing that $p$ has $\lambda$-DOP.
Let $q$ be any stationary, regular type non-orthogonal to $p$.
As $q$ is non-orthogonal to $M_3$, there is $q'\in S(M_3)$
non-orthogonal to $q$ (and hence to $p$) and conjugate to $q$.
But now, $q'\perp M_1$ and $q'\perp M_2$, so $(M_0,M_1,M_2,M_3)$
witnesses that $q'$ has $\lambda$-DOP.  Thus, $q$ has a DOP witness
by (2).

 (4) It suffices to prove this for $q$ lying directly above $p$.
As both notions are parallelism invariant, we may assume that
$p\in S(N)$, where $N$ is $\aleph_\eps$-saturated.   Choose $c$
realizing $p$
and $N[c]$ $\aleph_\eps$-prime over $N\cup\{c\}$ such that
$q\not\perp N[c]$, but $q\perp N$.   Choose $q'\in S(N[c])$ nonorthogonal 
to $q$.
Fix a cardinal $\lambda>|N|$, and choose a $\lambda$-DOP witness
$(M_0,M_1,M_2,M_3)$ for $p$.  Without loss, we may assume  that
$N\preceq M_3$ and that $\fg c N {M_3}$.   Let $M^*$ be $\lambda$-prime
over $N[c]\cup M_3$ and let
$q^*$ be the non-forking extension of $q'$ to $M^*$.  We argue that $(M_0,M_1,M_2,M^*)$
is a $\lambda$-DOP witness for $q^*$.  

To see this, first note that $N[c]$ is $\aleph_\eps$-constructible over
$N\cup\{c\}$, $N$ is $\aleph_\eps$-saturated, and $\fg c N {M_3}$,
so $N[c]$ is $\aleph_\eps$-constructible (hence $\lambda$-constructible)
over $M_3\cup\{c\}$.  Since $M^*$ is $\lambda$-constructible over 
$N[c]\cup M_3$, it follows that $M^*$ is $\lambda$-constructible over 
$M_3\cup\{c\}$,
hence is $\lambda$-prime over $M_3\cup\{c\}$.  Thus, by
Corollary~\ref{isom}, $M^*$ is $\lambda$-prime over $M_1\cup M_2$.
That is, $(M_0,M_1,M_2,M^*)$ is a $\lambda$-quadruple.

As well,  $q'\in S(N[c])$ is orthogonal to $N$ and $\fg {N[c]} N {M_3}$,
so $q'\perp M_3$.  As $M_1\cup M_2\subseteq M_3$, it follows
immediately that $q^*\perp M_1$ and $q^*\perp M_2$.
\endproof

Throughout the remainder of this paper, we consider sets $\P$ of stationary,
regular types over small subsets of the monster model $\C$.
We typically require $\P$ to be closed under automorphisms of $\C$ and nonorthogonality.

\begin{Definition}{\em
Let $\Treg$ denote the set of all
stationary, regular types over small subsets of $\C$
and fix a subset $\P\subseteq\Treg$ that is closed under automorphisms of $\C$
and nonorthogonality.

As notation, 
\begin{itemize}
\item A stationary type $q$ is orthogonal to $\P$, written $q\perp\P$,
if $q$ is orthogonal to every $p\in\P$.  
 $\P^{\perp}=\{q\in\Treg:q\perp\P\}$;
\item  $\Pactive$ is the closure of $\P$ in $\Treg$ under automorphisms,
nonorthogonality, and supporting (i.e., if $p\in\Treg$ supports some $q\in\P$,
then $p\in\Pactive$;
\item $\Pdull=\Treg\setminus \Pactive$.
\end{itemize}
}
\end{Definition}

\begin{Definition}{\em  Let $\P\subseteq \Treg$ be any set of regular
types.  A theory $T$ has
{\em $\P$-NDOP} if no $p\in\P$ has a DOP witness.
}
\end{Definition}

The following Corollary is merely a restatement of Proposition~\ref{DOPclosure}.

\begin{Corollary}  \label{DOPtransfer}  For any $\P\subseteq\Treg$,
$T$ has $\P$-NDOP if and only if $T$ has $\Pactive$-NDOP.
\end{Corollary}

\begin{Definition}  \label{Pdepth}  
{\em 
Given a class $\P$ of regular types, we define the {\em $\P$-depth\/} of a 
stationary, regular type $p$, $\dpP(p)\in{\bf ON}\cup\{-1\}$,  by
(1) $\dpP(p)=-1$ if and only if $p\in\Pdull$; and (2) $\dpP(p)\ge\alpha$
if and only if $p\in\Pactive$ and for every $\beta\in\alpha$ there is a triple
$(M,N,a)$, where $M$ is $\aleph_\eps$-saturated, $N$ is $\aleph_\eps$-prime over
$M\cup\{a\}$, $p$ is parallel to $\tp(a/M)$, and there is $q\in S(N)$ orthogonal to $M$
with $\dpP(q)\ge\beta$.
}
\end{Definition}

As in Chapter~X of \cite{Shc}, in the preceding definition it would be equivalent to
replace `$\aleph_\eps$-saturation' by `$\lambda$-saturation' for any uncountable cardinal $\lambda$.
The proof of the following Lemma is identical to the proof of Lemma~X~7.2 of \cite{Shc}.

\begin{Lemma}  \label{positivedepth}
If $T$ has $\P-NDOP$, then any regular $p$ with $\dpP(p)>0$ is trivial,
i.e., the set $p(\C)$ has a trivial pre-geometry with respect to the dependence relation of
forking.
\end{Lemma}

We close this section with two technical Lemmas that will be used later.
Note that a type $q$ (not necessarily regular) is orthogonal to $\Pdull$
if and only if every regular type non-orthogonal to $q$ is an element of $\Pactive$.

\begin{Lemma}  \label{oneone}
($\P$-NDOP, $\lambda\ge\aleph_\eps$) 
Suppose that $M$ is $\lambda$-prime over an independent triple
$(M_0,M_1,M_2)$ of $\lambda$-saturated models, $a$ is
$\eps$-finite satisfying $\tp(a/M)\perp\Pdull$ and $\tp(a/M)\perp M_2$.
Let $M[a]$ be any $\lambda$-prime model over $M\cup\{a\}$.
For any subset $N\subseteq M[a]$ that is maximal such that $\fg N {M_1} M$
we have:
\begin{enumerate}
\item  $N\preceq M[a]$, $N$ is $\lambda$-saturated, and $M[a]$ is
$\lambda$-prime over $N\cup M$; and
\item  For any $a^*\subseteq N$ such that $\fg N {M_1a^*} a$,
$N$ is $\lambda$-prime over $M_1\cup\{a^*\}$.
\end{enumerate}
\end{Lemma}

\bp  To see that $N\preceq M[a]$ and $N$ is $\lambda$-saturated,
choose $N^+\preceq M[a]$ to be $\lambda$-prime over $N$.
As $M_1$ is $\lambda$-saturated, it follows from Fact~\ref{refer}(2)
that $N^+$ is dominated by $N$ over $M_1$, hence $\fg {N^+} {M_1} M$,
so $N^+=N$ by the maximality of $N$.

Next, choose $M^*\preceq M[a]$ to be maximal such that $M^*$ is 
$\lambda$-saturated and $\lambda$-atomic over $N\cup M$.
(Since $T$ is superstable, the union of a continuous chain of 
$\lambda$-saturated models is $\lambda$-saturated, so $M^*$ exists.)
Since $a$ is $\eps$-finite, any subset $I\subseteq M[a]$
that is indiscernible over $M$ has
size at most $\lambda$ (when $\lambda=\aleph_\eps$, $I$ must be countable).
 It follows at once that every subset $I\subseteq M^*$
that is indiscernible over $N\cup M$ has size at most $\lambda$,
so by Fact~\ref{refer}(1) $M^*$ is $\lambda$-prime over $N\cup M$.
We complete the proof of (1) by showing that $M^*=M[a]$.

Suppose not.  Choose $c\in M[a]\setminus M^*$ such that $q=\tp(c/M^*)$
is regular.  The argument splits into cases.  First, if $q\perp N$ and $q\perp M$,
then $(M_1,N,M,M^*)$ is a DOP witness for $q$, so by
Corollary~\ref{isom}, any $\lambda$-prime model over $M^*\cup\{c\}$
is $\lambda$-prime over $N\cup M$, which contradicts the maximality
of $M^*$.  Second, if $q\not\perp N$, then choose a regular $r\in S(M^*)$
that does not fork over $N$ but $q\not\perp r$.  Choose $d\in M[a]\setminus M^*$
realizing $r$.  Then, by symmetry and transitivity of non-forking,
$\fg {Nd} {M_1} M$, which contradicts the maximality of $N$.
Finally, suppose that $q\not\perp M$.  As before, there is a regular $p\in S(M^*)$
that does not fork over $M$ but $q\not\perp p$, and an element $e\in M[a]\setminus M^*$ realizing $p$.   As $p$  is regular, based on $M$,
and non-orthogonal to $\tp(a/M)$, $p\in\Pactive$ and $p\perp M_2$. 
So, by $\P$-NDOP it must be that $p\not\perp M_1$.
But then, $p\not\perp N$, so arguing as above we contradict the maximality of
$N$.  This proves (1).

For (2), choose any such $a^*$.  We show that $N$ is $\lambda$-prime
over $M_1\cup\{a^*\}$ via Fact~\ref{refer}(1).  We already know
that $N$ is $\lambda$-saturated.  To see that $N$ is $\lambda$-atomic over
$M_1\cup\{a^*\}$, choose any finite set $c$ from $N$.  
As $N\subseteq M[a]$, $\tp(c/Ma)$ is $\lambda$-isolated.
But $\fg c {M_1a^*} {Ma}$, so $\tp(c/M_1a^*)$ is $\lambda$-isolated
as well (see e.g., \cite{Shc}~IV~4.1).
Finally, if $I\subseteq N$ is indiscernible over $M_1\cup\{a^*\}$,
then $I$ is indiscernible over $M_1$.  But $\fg N {M_1} M$, so
$I$ is indiscernible over $N\cup M$.  As $M[a]$ is $\lambda$-prime over
$N\cup M$, it follows that $I$ has size at most $\lambda$, completing
the proof of (2).
\endproof

\begin{Lemma} \label{twotwo}
($\P$-NDOP, $\lambda\ge\aleph_\eps$)  Suppose that $M_1\preceq M$
are both $\lambda$-saturated, $a$ is $\eps$-finite, $\tp(a/M)\perp\Pdull$, and 
either $\tp(a/M)$ does not fork over $M_1$, or $\tp(a/M)$ is regular
and non-orthogonal to $M_1$.  Let $M[a]$ be any $\lambda$-prime
model over $M\cup\{a\}$.
For any subset $N\subseteq M[a]$ that is maximal such that $\fg N {M_1} M$
we have:
\begin{enumerate}
\item  $N\preceq M[a]$, $N$ is $\lambda$-saturated, and $M[a]$ is
$\lambda$-prime over $N\cup M$; and
\item  For any $a^*\subseteq N$ such that $\fg N {M_1a^*} a$,
$N$ is $\lambda$-prime over $M_1\cup\{a^*\}$.
\end{enumerate}
\end{Lemma}

\bp  The proof is similar to the proof of Lemma~\ref{oneone}, only easier.
The hypotheses on $\tp(a/M)$ ensure that  for any $e\in M[a]\setminus M$,
as $e$ is dominated by $a$ over $M$, it follows that $\tp(e/M)\not\perp M_1$.

To see (1), take $N^+\preceq M[a]$ to be $\lambda$-prime over
$N$.  As before, the maximality of $N$ implies that $N^+=N$, so
$N\preceq M[a]$ and $N$ is $\lambda$-saturated.  As well, choose 
$M^*\preceq M[a]$ that is maximal such that $M^*$ is $\lambda$-saturated
and $\lambda$-atomic over $N\cup M$.  As before, indiscernible subsets
of $M^*$ over $N\cup M$ have size at most $\lambda$, so $M^*$
is $\lambda$-prime over $N\cup M$.  

The verification that $M^*=M[a]$ is also similar.  If not, choose 
$c\in M[a]\setminus M^*$ such that $q=\tp(c/M^*)$ is regular.
If $q\perp N$ and $q\perp M$, then $(M_1,N,M,M^*)$ is a DOP witness for $q$,
which again contradicts the maximality of $M^*$ by Corollary~\ref{isom}.
If $q\not\perp N$, then arguing as before there is a regular $r\in S(M^*)$
that does not fork over $N$, $q\not\perp r$, and a realization $d$ of $r$,
which contradicts the maximality of $N$.  Finally, if $q\not\perp M$,
then there is a regular $p\in S(M^*)$ that does not fork over $M$ but
$q\not\perp p$ and a realization $e$ of $p$ in $M[a]$.  Our conditions on
$\tp(a/M)$ imply that $\tp(e/M)\not\perp M_1$, hence $\tp(e/M)\not\perp N$
and we argue as above, completing the verification of (1).
The verification of (2) is identical to its verification in the proof of 
Lemma~\ref{oneone}.
\endproof

\section{$\P$-decompositions}
Throughout this section, assume that $T$ is superstable, and that $\P$ is a class of
regular types, closed under automorphisms of $\C$ and non-orthogonality.
We  define a number of species of $\P$-decompositions, along with a number of ways
in which one $\P$-decomposition can extend another.

\begin{Definition}  {\em
 Fix a model $M$.  A {\em weak $\P$-decomposition inside $M$}
is a sequence $\d=\<N_\eta,a_\eta:\eta\in I\>$ indexed by a tree 
$(I,\trianglelefteq)$
satisfying:
\begin{enumerate}
\item  $\{N_\eta:\eta\in I\}$ is an independent tree of elementary submodels of $M$;
\item  $\eta\trianglelefteq\nu$ implies $N_\eta\preceq N_\nu$;
\item Each $a_\eta\in N_\eta$ (but $a_{\<\>}$ is meaningless);
\item  For all $\nu\in Succ_I(\eta)$, $N_\nu$ is dominated by $a_\nu$ over $N_\eta$;
\item  If $\eta\neq\<\>$, then $\tp(a_\nu/N_\eta)\perp N_{\eta^-}$ for each $\nu\in Succ_I(\eta)$;
\item  For each $\eta\in I$, $\{a_\nu:\nu\in Succ_I(\eta)\}$ is independent over $N_\eta$
and $\tp(a_\nu/N_\eta)\perp \Pperp$ for each $\nu\in Succ_I(\eta)$.
\end{enumerate}
}
\end{Definition}

Note that in the Definition above, we do not require that $\tp(a_\nu/N_\eta)$ be
regular.  However, the content of (6) is that
any regular type $q\not\perp\tp(a_\nu/N_\eta)$ is necessarily in $\P$.

\begin{Lemma}  \label{easy}
Suppose $\d=\<N_\eta,a_\eta:\eta\in I\>$ is a weak $\P$-decomposition inside $M$.
Then:
\begin{enumerate}
\item  If $I_1,I_2\subseteq I$ are both downward closed and $I_0=I_1\cap I_2$, then
$$ \fg {\left(\bigcup_{\eta\in I_1} N_\eta\right)}  
{\left(\bigcup_{\eta\in I_0} N_\eta\right)}  
{\left(\bigcup_{\eta\in I_2} N_\eta\right)}  
$$
\item  If $\eta\in I$, $\nu=\eta\conc\<\alpha\>$, where $\alpha$ is least such that $\eta\conc\<\alpha\>\not\in I$,
the element $a_\nu\in M$ satisfies $\tp(a_\nu/N_\eta)\perp\Pperp$, if $\eta\neq\<\>$ then $\tp(a_\nu/N_\eta)\perp N_{\eta^-}$,
and $\fg {a_\nu}  {N_\eta} {\{a_\gamma:\gamma\in Succ_I(\eta)\}}$, and $N_\nu\preceq M$ is dominated by $a_\nu$ over $N_\eta$,
then $\d^*=\d\conc\<N_\nu,a_\nu\>$ is a weak $\P$-decomposition inside $M$.
\end{enumerate}
\end{Lemma}

There are two ways of defining when a weak $\P$-decomposition inside
a model $M$ is `maximal'.    Fortunately, at least when both $M$ and each of the submodels $N_\eta$ are $\aleph_\eps$-saturated, Lemma~\ref{stopping}
shows that the two notions are equivalent.  

\begin{Definition} \label{maximaldecomp} {\em  Suppose that
  $\d=\<N_\eta,a_\eta:\eta\in I\>$  of $M$ is a weak
$\P$-decomposition inside $M$.  As notation, for each $\eta\in I$, let
$$C_\eta(M)=\{a\in M\setminus N_\eta:\tp(a/N_\eta)\perp\Pperp\ 
\hbox{and $\perp N_{\eta^-}$ (when $\eta\neq\<\>$)}\}$$
\begin{itemize}
\item  {\em $\d$ is a weak $\P$-decomposition {\bf  of ${\bf M}$}} if, for every $\eta\in I$, $\{a_\nu:\nu\in Succ_I(\eta)\}$ is a maximal $N_\eta$-independent subset of
$C_\eta(M)$.
\item  {\em $\d$ $\P$-exhausts $M$} if, for every $\eta\in I$  for every
regular $p\in S(N_\eta)\cap\P$ orthogonal to $N_{\eta^-}$ (when $\eta\neq\<\>$)
and for every $e\in p(\C)$, if $\fg e {N_\eta} {\{a_\nu:\nu\in Succ_I(\eta)\}}$
then $\fg e {N_\eta} M$.
\end{itemize}
}
\end{Definition}

\begin{Lemma}  \label{stopping}
Suppose that $\d=\<N_\eta,a_\eta:\eta\in I\>$ is a 
weak $\P$-decomposition inside an $\aleph_\eps$-saturated model $M$
such that every $N_\eta$ is $\aleph_\eps$-saturated as well.
Then $\d$ is weak $\P$-decomposition of $M$ if and only if
$\d$ $\P$-exhausts $M$.
\end{Lemma}

\bp  For both directions, recall that if $h\in M\setminus N_\eta$, then there is
a finite, $N_\eta$-independent set $\{b_i:i<n\}\subseteq M$ domination equivalent to $h$ over $N_\eta$ with
$\tp(b_i/N_\eta)$
is regular for each $i<n$.  

For the easy direction, suppose that $\d$ is not a weak $\P$-decomposition of $M$.
Choose $\eta\in I$ such that $A=\{a_\nu:\nu\in Succ_I(\eta)\}$ is not maximal in
$C_\eta(M)$.  Choose $h\in C_\eta(M)$ such that $\fg h {N_\eta} A$, and from
above, choose $\{b_i:i<n\}\subseteq M$ domination equivalent to $h$ over 
$N_\eta$.  Then, for any $i<n$, the element $b_i$ and the type $\tp(b_i/N_\eta)$ witness that $\d$ does not $\P$-exhaust $M$.

Conversely, suppose that $\d$ is a weak $\P$-decomposition of $M$.  Fix any $\eta\in I$, any regular type $p\in S(N_\eta)\cap\P$ that is $\perp N_{\eta^-}$
when $\eta\neq\<\>$, and any $e\in p(\C)$ with $\nfg e {N_\eta} M$.
We will show that $\nfg e {N_\eta} A$, where $A=\{a_\nu:\nu\in Succ_I(\eta)\}$.

To see this, using the note above choose $n<\omega$ minimal such 
that there are $h\in M\setminus N_\eta$ and $B=\{b_i:i<n\}\subseteq M$
such that $\nfg e {N_\eta} h$, and $h$ and $B$ are domination equivalent over $N_\eta$ with $\tp(b_i/N_\eta)$ regular for each $i$.  It follows
from the minimality of $n$ that $\tp(b_i/N_\eta)$ is non-orthogonal to $p$,
hence each $b_i\in C_\eta(M)$.  As $A$ is maximal in $C_\eta(M)$, we have
that $\nfg {b_i} {N_\eta} A$ for each $i$, hence $\tp(b_i/N_\eta A)$ is hereditarily orthogonal to $p$ 
(i.e., $\tp(b_i/N_\eta A)$, as well as every forking extension of it is orthogonal to $p$).  Thus, $\tp(B/N_\eta A)$ is hereditarily orthogonal to $p$.  
This implies $\nfg e {N_\eta} A$.  [Why?  If not, then $\tp(e/N_\eta A)$ would be parallel to $p$, so by orthogonality we would have $\fg e {N_\eta A} B$.
This would imply that $e$ and $B$ (and hence $e$ and $h$) are independent over $N_\eta$, which is a contradiction.]
\endproof
\endproof

For our next series of results, we insist that the model $M$ be sufficiently saturated, and
we additionally require that each submodel occurring in a decomposition be sufficiently saturated as well.
In most applications, $\aleph_\eps$-saturation would suffice, but it costs little to work
in the more general context of $\lbarP$-saturated models, which we now introduce.

\medskip
{\bf Fix, for the remainder of this section, a pair ${\bf \lbar=(\lambda,\mu)}$ of cardinals satisfying
${\bf \lambda,\mu\ge\aleph_\epsilon}$.}  Throughout the whole of this paper, if $\lambda=\mu=\aleph_\eps$, we write
$\aeP$ in place of $(\aleph_\eps,\aleph_\eps,\P)$.

\medskip

\begin{Definition}
{\em
We say that a model $M$ is {\em $\lbarP$-saturated} if it is $\aleph_\eps$-saturated, 
and for each finite $A\subseteq M$, $\dim(p,M)\ge\lambda$
for each $p\in\P\cap S(A)$, and $\dim(q,M)\ge\mu$ for all
stationary, regular $q\in\Pperp\cap S(A)$.  (If either $\lambda$ or $\mu$ is $\aleph_\eps$, the associated dimension is at least $\aleph_0$.)

We say that a $\lbarP$-saturated  model $N$ is {\em $\lbarP$-prime over a set $X$}
if $N\contains X$ and 
$N$ embeds elementarily over $X$ into any $\lbarP$-saturated model containing $X$.  
}
\end{Definition}

Note that our assumptions on $\lbar$ guarantee that any $\lbarP$-saturated model is $\aleph_\eps$-saturated, but we include this clause for emphasis.
Also, if $\lambda=\mu$, then the $\lbarP$-saturated models are precisely the
$\lambda$-saturated models.   
The standard facts about the existence $\lbarP$-prime models extend easily to this context.  
To see this, call a type $\tp(a/B)$ {\em $\lbarP$-isolated} if any of the three conditions hold: (1)
$\tp(a/B)$ is $\aleph_\eps$-isolated ($={\bf F}^a_{\aleph_0}$-isolated)  or (2) $\tp(a/B)\in \P$ and
is $\lambda$-isolated; or (3) $\tp(a/B)\in\Pperp$ and is $\mu$-isolated.  Next, call a set $B$ {\em $\lbarP$-primitive over $A$}
if 
$B=A\cup\{b_i:i<\alpha\}$, where $\tp(b_i/A\cup\{b_j:j<i\})$ is $\lbarP$-isolated
for every $i$, and call a model $M$ {\em $\lbarP$-primary over $A$} if $M$ is $\lbarP$-saturated and its universe is
$\lbarP$-primitive over $A$.  This notion of isolation satisfies the same axioms as
for ${\bf F}^a_\lambda$-isolation in Chapter~4 of \cite{Shc} and thus
we obtain the same consequences.  In particular:
\begin{itemize}
\item  If $A\subseteq M^*$ with $M^*$ $\lbarP$-saturated, then there is a $M\preceq M^*$ that is $\lbarP$-primary over $A$;
\item  $M$
$\lbarP$-primary over $A$ implies $M$ is $\lbarP$-prime over $A$; and
\item If $M$ is $\lbarP$-saturated, $M\subseteq A$, and $N$ is $\lbarP$-prime over $A$, then $N$ is dominated by $A$ over $M$.
\end{itemize}

\begin{Definition}  {\em 
Suppose that $M$ is $\lbarP$-saturated.  A {\em weak $\lbarP$-decomposition inside $M$ (of $M$)}
is a weak $\P$-decomposition inside $M$ (of $M$) for which each of the submodels $N_\eta$ is an $\lbarP$-saturated
elementary substructure of $M$.   
} 
\end{Definition}

A salient feature of weak $\lbarP$-decompositions is that each of the
submodels is itself  $\aleph_\epsilon$-saturated.
The proof of the following Lemma is virtually identical to arguments
in Section~X.3 of \cite{Shc}.

\begin{Lemma}[$\P$-NDOP]  \label{X.3}
Suppose that $\<N_\eta,a_\eta:\eta\in I\>$
is a weak $\lbarP$-decomposition of a $\lbarP$-saturated model
$M$.  Let $\Mbar\preceq M$ be
any $\aleph_\epsilon$-prime submodel of $M$ over 
$\bigcup_{\eta\in I} N_\eta$.
Then
if $p\in\P$ is non-orthogonal to $\Mbar$, then there is a unique,
$\triangleleft$-minimal $\eta\in I$ such that $p\not\perp N_\eta$.
\end{Lemma}

\bp  We first show that
$p\not\perp N_\eta$ for some $\eta\in I$.
As $\Mbar$ is $\aleph_\epsilon$-saturated, there is $q\in S(\Mbar)$
that is regular and non-orthogonal to $p$.  As any such $q$ is in $\P$,
we may assume that $p\in S(\Mbar)$ to begin with.
Choose a finite $B\subseteq\Mbar$ over which $p$ is based and stationary.
As $B$ is $\aleph_\epsilon$-isolated over $\bigcup_{\eta\in I} N_\eta$,
there is a finite subtree $I_0\subseteq I$ such that $B$ is $\aleph_\epsilon$-isolated
over $\bigcup_{\eta\in I_0} N_\eta$.  Choose any $M_0\preceq \Mbar$
such that $B\subseteq M_0$ and $M_0$ is $\aleph_\epsilon$-prime over
$\bigcup_{\eta\in I_0} N_\eta$.  As there is some type $p'\in S(M_0)$
parallel to $p$,
 $\P$-NDOP implies that $p\not\perp N_\eta$ for some $\eta\in I_0$.

Finally,  using Lemma~\ref{easy}(1) it follows that there is a unique 
$\triangleleft$-minimal $\eta\in I$ with $p\not\perp N_\eta$.
\endproof

The following definition makes sense in our context, as $\lbarP$-decompositions
 have no
control over types orthogonal to $\P$.

\begin{Definition} {\em  An $\ae$-saturated model $N$ is {\em $\P$-minimal over $X$}
if $N\contains X$, but for any $\ae$-saturated $N_0\preceq N$ containing $X$,
$\tp(e/N_0)\perp \P$ for every $e\in N\setminus N_0$.
}
\end{Definition}

\begin{Corollary}[P-NDOP]  \label{Pminimal}
Suppose that $\<N_\eta,a_\eta:\eta\in I\>$
is a weak $\lbarP$-decomposition of a $\lbarP$-saturated model
$M$ and let $\Mbar\preceq M$ be any
$\ae$-prime model over
$\bigcup_{\eta\in I} N_\eta$.
Then:
\begin{enumerate}
\item Every $c\in M\setminus \Mbar$ satisfies $\tp(c/\Mbar)\perp\P$; and
\item  $\Mbar$ is  $\P$-minimal over 
$\bigcup_{\eta\in I} N_\eta$.
\end{enumerate}
\end{Corollary}

\bp  
(1) Assume by way of contradiction that there is $c\in M$ such that
$\tp(c/\Mbar)\not\perp\P$.  As $\P$ is closed under non-orthogonality and
automorphisms of $\C$, there is $p\in\P\cap S(\Mbar)$ non-orthogonal to
$\tp(c/\Mbar)$.  Then, by Fact~\ref{refer}(3), there is $e\in M$
realizing $p$.
So, by Lemma~\ref{X.3}, $p\not\perp N_\eta$ for some $\eta\in I$.
Thus, by Fact~\ref{refer}(4) there is $e^*\in M$
domination equivalent to $e$ over $\Mbar$ with $\fg {e^*} {N_\eta} {\Mbar}$.
As $\{a_\nu:\nu\in Succ_I(\eta)\}\subseteq \Mbar$, this contradicts
the fact that $\<N_\eta,a_\eta:\eta\in I\>$
is a weak $\lbarP$-decomposition of $M$.

(2)  Choose any $M_1\preceq\Mbar$ that is $\ae$-prime over 
$\bigcup_{\eta\in I} N_\eta$.  Then (1) applies to $M_1$.  That is,
there is no $c\in M\setminus M_1$ such that $\tp(c/M_1)\not\perp\P$.
Thus, $\Mbar$ is $\P$-minimal over $\bigcup_{\eta\in I} N_\eta$.
\endproof

Next, we show that if we additionally assume that $\P=\Pactive$,
then we can extend the previous results to any $\ae$-saturated submodel
of $M$ containing the decomposition.

\begin{Proposition}[$\P$-NDOP, $\P=\Pactive$]  \label{noPrealized}
Suppose that $\<N_\eta,a_\eta:\eta\in I\>$
is a weak $\lbarP$-decomposition of a $\lbarP$-saturated model $M$.  
Let $M^*\preceq M$ be
any $\aleph_\epsilon$-saturated model containing
$\bigcup_{\eta\in I} N_\eta$.
Then there is no $e\in M\setminus M^*$ such that $\tp(e/M^*)\not\perp\P$.
\end{Proposition}

\bp  As both $M^*$ and $M$ are $\aleph_\epsilon$-saturated, it suffices
to prove that there is no $e\in M\setminus M^*$ such that $\tp(e/M^*)\in\P$.
Assume by way of contradiction that there were such an $e$.
Let $M_0\preceq M^*$ be any $\aleph_\epsilon$-prime model over
$\bigcup_{\eta\in I} N_\eta$.  Next, form an increasing
sequence $\<M_\alpha:\alpha\le\delta\>$ of $\aleph_\epsilon$-saturated
models, with $M_\delta=M^*$, $M_{\alpha+1}$ is $\aleph_\epsilon$-prime
over $M_\alpha\cup\{b_\alpha\}$, where $\tp(b_\alpha/M_\alpha)$ is regular,
and for $\alpha<\delta$ a non-zero limit, $M_\alpha$ is 
$\aleph_\epsilon$-prime over $\bigcup_{\beta<\alpha} M_\beta$.

Choose $\alpha\le\delta$ least such that there is some $e\in M\setminus M_\alpha$
such that $\tp(e/M_\alpha)\in\P$.
By superstability, $\alpha$ cannot be a non-zero limit ordinal.
Now suppose $\alpha=\beta+1$.   On one hand, if  $p=\tp(e/M_\alpha)\in\P$ were non-orthogonal to $M_\beta$,
then by Fact~\ref{refer}(4), there would be $e^*\in M$ such that $q=\tp(e^*/M_\beta)$ is regular and non-orthogonal
to $p$, contradicting the minimality of $\alpha$.  On the other hand, if $p\perp M_\beta$,
then
as $\Pactive=\P$, $r=\tp(b_\beta/M_\beta)\in \P$,  which again contradicts the minimality of $\alpha$.

Thus, $\alpha$ must equal zero, i.e., there is $e\in M\setminus M_0$
such that $p=\tp(e/M_0)\in\P$.  By Lemma~\ref{X.3},
choose a $\triangleleft$-minimal $\eta\in I$ such that $p\not\perp N_\eta$.

 Choose $q\in S(N_\eta)$ regular such that $p\not\perp q$
and let $q'\in S(M_0)$ be the non-forking extension of $q$ to $M_0$.
As both $M_0$ and $M$ are $\aleph_\epsilon$-saturated, there is
$c\in M\setminus M_0$ realizing $q'$.  As $q'\in\P$, we have $c\in C_\eta(M)$
in the notation of Definition~\ref{maximaldecomp}, which contradicts the maximality
of $\{a_\nu:\nu\in Succ_I(\eta)\}$.
\endproof

\begin{Corollary}[$\P$-NDOP, $\P=\Pactive$]  \label{nonorthtransfer}
Suppose that $\<N_\eta,a_\eta:\eta\in I\>$
is a weak $\lbarP$-decomposition of $M$.  Let $M^*\preceq M$ be
any $\aleph_\epsilon$-saturated elementary submodel containing
$\bigcup_{\eta\in I} N_\eta$.
If $p\in\P$ and $p\not\perp M$, then $p\not\perp M^*$.
\end{Corollary}  

\bp  
As in the proof above, form an increasing
sequence $\<M_\alpha:\alpha\le\delta\>$ of $\aleph_\epsilon$-saturated
models,  this time with $M_0=M^*$, $M_\delta=M$, $M_{\alpha+1}$ is $\aleph_\epsilon$-prime
over $M_\alpha\cup\{b_\alpha\}$, where $\tp(b_\alpha/M_\alpha)$ is regular,
and for $\alpha<\delta$ a non-zero limit, $M_\alpha$ is 
$\aleph_\epsilon$-prime over $\bigcup_{\beta<\alpha} M_\beta$.
Choose $\alpha\le\delta$ least such that $p\not\perp M_\alpha$.
We will show that $\alpha=0$.  Clearly, $\alpha$ cannot be a 
non-zero limit by superstability.  Assume by way of contradiction
that $\alpha=\beta+1$.  Then $p\not\perp M_\alpha$, but $p\perp M_\beta$.
But, as before, this implies that $r=\tp(b_\beta/M_\beta)\in \Pactive=\P$.
But now, $M_\beta$ is an $\aleph_\epsilon$-saturated model containing
$\bigcup_{\eta\in I} N_\eta$, yet there is an element of $M\setminus M_\beta$
realizing $r\in\P$, contradicting Proposition~\ref{noPrealized}.
Thus, $\alpha=0$, so $p\not\perp M^*$.
\endproof

\begin{Corollary}[$\P$-NDOP, $\P=\Pactive$]  \label{mineta}
Suppose that $\<N_\eta,a_\eta:\eta\in I\>$
is a weak $\lbarP$-decomposition of a $\lbarP$-saturated model
$M$.  
If $p\in\P$ and $p\not\perp M$, then 
there is a unique $\triangleleft$-minimal $\eta\in I$ such that $p\not\perp N_\eta$.
\end{Corollary}  

\bp  Let $M^*\preceq M$ be any $\aleph_\epsilon$-prime model over
$\bigcup_{\eta\in I} N_\eta$.  By Corollary~\ref{nonorthtransfer}
$p\not\perp M^*$, so by Lemma~\ref{X.3}, $p\not\perp N_\eta$ for some
$\triangleleft$-minimal $\eta\in I$.
As in the proof of Lemma~\ref{X.3}, the uniqueness follows from 
Lemma~\ref{easy}(1).
\endproof

Until this point in our discussion, the submodels occurring in a decomposition
could be very large, with an extreme case being that any model $M$
has a one-element decomposition $\<M\>$.  The next definition limits the size
of the submodels, while retaining the fact that they are at least $\aleph_\eps$-saturated.

\begin{Definition}  {\em
A {\em prime $\lbarP$-decomposition inside $M$ (of $M$)}
is a weak $\lbarP$-decomposition inside $M$ (of $M$) in which
$N_{\<\>}$ is $\lbarP$-prime over $\emptyset$ and, for each 
$\eta\in I\setminus \{\<\>\}$, $N_\eta$ is $\lbarP$-prime over
$N_{\eta^-}\cup\{a_\eta\}$.
}
\end{Definition}

\begin{Definition}  {\em  Fix a $\lbarP$-saturated model $M$.  A prime 
$\lbarP$-decomposition $\d_2=\<N_\eta^2,a_\eta^2:\eta\in J\>$
{\em end extends} the prime $\lbarP$-decomposition
$\d_1=\<N_\eta^1,a_\eta^1:\eta\in I\>$ if $I\subseteq J$ and,
for each $\eta\in I$, $N^2_\eta=N^1_\eta$ and $a_\eta^2=a_\eta^1$.

We say $\d_2$ is a {\em regular end extension} of $\d_1$ if, in addition
$\tp(a_\eta/N_{\eta^-})$ is regular for each $\eta\in J\setminus I$.
Furthermore, $\d_2$ is a {\em standardly regular end extension} of $\d_1$
if, $\tp(a_\eta/N_{\eta^-})=\tp(a_\nu/N_{\nu^-})$
whenever  $\eta,\nu$ in $J\setminus I$, $\eta^-=\nu^-$, and
$\tp(a_\eta/N_{\eta^-})\not\perp\tp(a_\nu/N_{\nu^-})$.
}
\end{Definition}

The following Lemma is straightforward, and relies on the fact that
if $N\preceq M$ are both $\lbarP$-saturated with $a\in M\setminus N$ 
satisfying $\tp(a/N)\in\P$, then there is $N[a]\preceq M$ that is
$\lbarP$-prime over $N\cup\{a\}$ and that $N[a]$ contains realizations
of every regular type over $N$ non-orthogonal to $\tp(a/N)$.
Proofs of similar statements appear in Section~X.3 of \cite{Shc}.

\begin{Lemma}  \label{maxendextension}
Suppose $\d=\<N_\eta,a_\eta:\eta\in I\>$ is a prime $\lbarP$-decomposition
inside an $\lbarP$-saturated model $M$.  Then:
\begin{enumerate}
\item $\d$ is a prime $\lbarP$-decomposition of $M$ if and only if it
has no proper (standardly regular) end extension; and
\item  There is a prime $\lbarP$-decomposition $\d^*$ of $M$ that is
a standardly regular end extension of $\d$.
\end{enumerate}
\end{Lemma}

%
%

Similarly to the main theme of \cite{Sh401}, we wish to investigate
$\P$-decompositions that lie above a specific triple
$(N,N',a)$, where $N\preceq N'$, with $\tp(a/N')\perp\Pperp$ and
$\tp(a/N')\perp N$.  That is, triples where $a$ could play the role of
$a_\nu$ in some $\P$-decomposition with $N'=N_\eta$ and
$N=N_{\eta^-}$.  However, as in \cite{Sh401}, this is too much
data to record at once, so we seek an $\epsilon$-finite approximation of
it.

Specifically, for $M$ any model,  let
$$\Gamma(M):=
\{(A,B):A\subseteq B\subseteq M\ \hbox{are both $\epsilon$-finite}\}$$
We frequently write $B\choose A$ for elements of $\Gamma(M)$, and
if $A$ is not a subset of $B$, we mean ${A\cup B}\choose A$.  Let

$$\Gamma_\P(M):=\{{B\choose A}\in\Gamma(M):\tp(B/A)\perp\Pperp\}$$

\begin{Definition} \label{ofMoverBA}
{\em  For ${B\choose A} \in \Gamma_\P(M)$, 
a {\em prime $\lbarP$-decomposition over $B\choose A$ inside $M$,}
$\d=\<N_\eta,a_\eta:\eta\in I\>$, is a prime $\lbarP$-decomposition inside $M$
in which $\<0\>$ is the unique successor of $\<\>$ in $I$,
$A\subseteq N_{\<\>}$, $B\subseteq N_{\<0\>}$, and 
$B\subseteq\dcl(a_{\<0\>})$.   By analogy with 
Definition~\ref{maximaldecomp},
\begin{itemize}
\item  such a {\em $\d$ is {\bf  of ${\bf M}$}} if, for every $\eta\in I\neq\<\>$, 
$\{a_\nu:\nu\in Succ_I(\eta)\}$ is a maximal $N_\eta$-independent subset of
$C_\eta(M)$; and
\item  {\em $\d$ $\P$-exhausts $M$ over ${B\choose A}$} if, for every $\eta\in I\neq\<\>$  for every
regular $p\in S(N_\eta)\cap\P$ orthogonal to $N_{\eta^-}$ (when $\eta\neq\<\>$)
and for every $e\in p(\C)$, if $\fg e {N_\eta} {\{a_\nu:\nu\in Succ_I(\eta)\}}$
then $\fg e {N_\eta} M$.
\end{itemize}
}
\end{Definition}

The following Lemma is straightforward.  The verification of (5) is analogous to the proof of Lemma~\ref{stopping}.

\begin{Lemma} \label{straightBA}
 Fix a $\lbarP$-saturated model $M$ and ${B\choose A} \in\Gamma_\P(M)$.
\begin{enumerate}  
\item  If $N_{\<\>}\preceq M$ is $\lbarP$-prime over $\emptyset$, contains
$A$, and $\fg B A {N_{\<\>}}$, and $N_{\<0\>}\preceq M$ is $\lbarP$-prime
over $N_{\<\>}\cup B$, then $\<N_{\<\>},N_{\<0\>}\>$ is
a prime $\lbarP$-decomposition over $B\choose A$ inside $M$;
\item  A prime $\lbarP$-decomposition $\d$ over $B\choose A$ inside $M$ is 
a prime $\lbarP$-decomposition over $B\choose A$ of $M$ if and only if
$\d$ has no proper $\lbarP$-decomposition over $B\choose A$ end extending it;

\item  Every prime $\lbarP$-decomposition over $B\choose A$ inside $M$
has a (standardly regular) end extension to  
a prime $\lbarP$-decomposition over $B\choose A$ of $M$;
\item  Every prime $\lbarP$-decomposition $\d$ over $B\choose A$ inside $M$
is a prime $\lbarP$-decomposition inside $M$, hence
has a (standardly regular) end extension to a 
prime $\lbarP$-decomposition $\d^*$ of $M$;
Moreover, if $\d$ is a decomposition over $B\choose A$ of $M$ and 
is indexed by the tree $(I,\triangleleft)$ and $\d^*$ is indexed
by $(J,\triangleleft)$, then $\neg(\<0\>\trianglelefteq \eta)$ for all $\eta\in J\setminus I$.
\item  A $\lbarP$-decomposition $\d$ over $B\choose A$ inside $M$ is of $M$ if and only if $\d$ $\P$-exhausts
$M$ over $B\choose A$.
\end{enumerate}
\end{Lemma}

\section{Trees of subsets of an $\aleph_\eps$-saturated model}

{\bf Throughout this section $T$ is superstable with $\P$-NDOP,
and $\P$ is closed under automorphisms of $\C$, non-orthogonality,
and $\P=\Pactive$.}  

In addition, all models $M$ we consider will be $\aleph_\eps$-saturated,
and all decompositions we consider will be $\aeP$-decompositions inside/of $M$.

\begin{Definition} \label{PP}
{\em
Fix an $\ae$-saturated
model $M$ and ${B\choose A}\in\Gamma_\P(M)$.
 Suppose $\d=\<N_\eta,a_\eta:\eta\in I\>$ is a prime
$\aeP$-decomposition over $B\choose A$ of $M$.
Then
$$\PP(\d,M)=\{p\in S(M):p\in\P, p\perp N_{\<\>},\ \hbox{but $p\not\perp
N_\eta$ for some $\eta\in I\setminus\{\<\>\}$}\}$$
}
\end{Definition}

The goal for this section will be Theorem~\ref{big}, which asserts
that $\PP(\d_1,M)=\PP(\d_2,M)$ for any two prime, $\aeP$-decompositions
$\d_1,\d_2$ of $M$ above $B\choose A$.
We begin by introducing
another way of `increasing' a decomposition.
\begin{Definition}  {\em
 A prime 
$\aeP$-decomposition $\d_2=\<N_\eta^2,a_\eta^2:\eta\in J\>$ inside $\C$
is a {\em blow up} of the prime $\aeP$-decomposition
$\d_1=\<N_\eta^1,a_\eta^1:\eta\in I\>$ inside $\C$ if $J=I$, but for every
$\eta\in I$, $N_\eta^1\preceq N_\eta^2$ and, when $\eta\neq\<\>$,
$N^2_\eta$ is $\aeP$-prime over $N^1_\eta\cup N^2_{\eta^-}$.
}
\end{Definition}


\begin{Lemma} \label{1.20}
Suppose that $M$ is $\ae$-saturated, ${B\choose A}\in \Gamma_\P(M)$, $\d_2=\<N^2_\eta,a_\eta:\eta\in I\>$ is a blow up of
$\d_1=\<N^1_\eta,a_\eta:\eta\in I\>$, $A\subseteq B\subseteq N^1_{\<0\>}$,
and each $N^2_\eta\preceq M$.
Then:
\begin{enumerate}
\item  If $\nu\in Succ_I(\eta)$, then $\fg {N^2_\eta} {N^1_\eta} {N_\nu^1}$;
\item  If $Y=\{\rho\in I:\neg(\eta\triangleleft\rho)\}$, then
$\fg {N^2_\eta} {N^1_\eta} {\bigcup_{\rho\in Y} N^1_\rho}$ and 
$\eta\triangleleft\nu$ implies
$\fg {N^2_\nu} {N^1_\eta} {\bigcup_{\rho\in Y} N^1_\rho}$;
\item $\d_2$ is an $\aeP$-decomposition inside $M$ above ${B\choose A}$ if
and only if $\d_1$ is; and
\item $\d_2$ is an $\aeP$-decomposition of $M$ above ${B\choose A}$ if
and only if $\d_1$ is.
\end{enumerate}
\end{Lemma}

\bp  This is exactly analogous to Fact~1.20 of \cite{Sh401}.
In the proof of (4), we need to appeal to $\P$-NDOP instead of NDOP.
\endproof

\begin{Lemma}  \label{blowupPP}
Suppose $M$ is $\ae$-saturated and 
${B\choose A}\in\Gamma_\P(M)$.
If $\d_2$ is a blow up of $\d_1$ and both $\d_1$, $\d_2$ are
$\aeP$-decompositions of $M$ over ${B\choose A}$, then
$\PP(\d_1,M)=\PP(\d_2,M)$.
\end{Lemma}

\bp  This is very much like Fact 1.22 of \cite{Sh401}, but we give details.
As notation, say $\d_\ell=\<N_\eta^\ell,a_\eta:\eta\in I\>$
for $\ell=1,2$.  Fix $p\in S(M)\cap \P$, so in particular, $p$ is regular.
We must prove that
$$\left(p\perp N^1_{\eta^-}\ \hbox{and}\ p\not\perp N^1_\eta\right)
\Leftrightarrow 
\left(p\perp N^2_{\eta^-}\ \hbox{and}\ p\not\perp N^2_\eta\right)$$
for every $\eta\neq\<\>$.

First, assume $\eta\neq\<\>$ and
$p\perp N^1_{\eta^-}$ and $p\not\perp N^1_\eta$.
As $N_\eta^1\preceq N_\eta^2$, $p\not\perp N^2_\eta$ trivially.
Also, choose  a regular $q\in S(N^1_\eta)$ with $p\not\perp q$.
Then $q\perp N^1_{\eta^-}$ since $p$ is, and it suffices to show that
$q$ is orthogonal to $N^2_{\eta^-}$.  But this follows immediately since
$\fg {N^1_\eta} {N^1_{\eta^-}} {N^2_{\eta^-}}$.

Conversely, assume
$\eta\neq\<\>$ and
$p\perp N^2_{\eta^-}$ and $p\not\perp N^2_\eta$.
Then, since $N^1_{\eta^-}\preceq N^2_{\eta^-}$, $p\perp N^1_{\eta^-}$.
As well, $(N^1_{\eta^-}, N^2_{\eta^-}, N^1_\eta)$ form an independent
triple of $\ae$-saturated models (see Definition~\ref{1.1})
and $N^2_\eta$ is $\ae$-prime over their union.
Thus, as $p\in\P$, it follows from $\P$-NDOP that $p\not\perp N^1_\eta$.
\endproof

\begin{Lemma}  \label{charPP}
Suppose that
$M$ is $\ae$-saturated, ${B\choose A}\in\Gamma_\P(M)$, and for $\ell=1,2$ $\d_\ell=\<N_\eta^\ell,a_\eta^\ell:\eta\in I_\ell\>$ are each prime, $\ae$-decompositions of $M$ above
$B\choose A$.  If $N^1_{\<\>}=N^2_{\<\>}$ then
$\PP(\d_1,M)=\PP(\d_2,M)$.
\end{Lemma}

\bp  First, by Lemma~\ref{straightBA}(4), choose a prime, $\ae$-prime decomposition
$\d^*_1=\<N_\eta^1,a_\eta^1:\eta\in J_1\>$
of $M$ end extending $\d_1$.  As notation, let $H=J_1\setminus I_1$
and for each $\eta\in H$, let $N^2_\eta=N^1_\eta$ and $a^2_\eta=a^1_\eta$.
It is easily checked that $\d_2^*:=\<N^2_\eta,a^2_\eta:\eta\in I_2\cup H\>$
is an $\aeP$-decomposition of $M$.  

Now, for each $p\in S(M)\cap\P$ with $p\perp N^1_{\<\>}$ and for each $\ell=1,2$ there is a unique
$\eta(p,\ell)\in I_\ell\cup H$ such that $p\not\perp N_{\eta(p,\ell)}$,
but $p\perp N_{\eta(p,\ell)^-}$.  But, as $N^2_\eta=N^1_\eta$ for each 
$\eta\in H$, $\eta(p,1)\in H$ if and only if $\eta(p,2)\in H$.

Thus, for each $p\in S(M)\cap \P$ that is orthogonal to $N^1_{\<\>}=N^2_{\<\>}$ we have
$p\in\PP(\d_1,M)$ if and only if $\eta(p,1)\in H$ if and only if $\eta(p,2)\in H$
if and only if $p\in\PP(\d_2,M)$.
\endproof

We come to the issue of the existence of blow-ups of decompositions.
It is comparatively easy to blow up a decomposition inside an $\ae$-saturated
model $M$.

\begin{Lemma}  \label{blowup}  Suppose that $M$ is $\ae$-saturated
and $\d=\<N_\eta,a_\eta:\eta\in I\>$ is a prime $\aeP$-decomposition
inside $M$.  For any $N^*$ satisfying $N_{\<\>}\preceq N^*\preceq M$
that is $\ae$-prime over $\emptyset$, there is a prime $\aeP$-decomposition
$\d^*$ inside $M$ with $N^{\d^*}_{\<\>}=N^*$ that is a blow up of $\d$.
\end{Lemma}

\bp  Choose any enumeration $\<\eta_i:i<i^*\>$ of $I$ such that
$\eta_i\triangleleft\eta_j$ implies $i<j$ and so that for some $\alpha^*\le i^*$
$\eta_i\in Succ_I(\<\>)$ if and only if $1\le i<\alpha^*$.
Note that $\eta_0=\<\>$ for any such enumeration.
Put $N^*_{\<\>}:=N^*$.  Then, by induction on $1\le i<i^*$,
argue that 
$$\fg {N^*} {N_{\<\>}} {\bigcup_{j<i} N^*_{\eta_j}}$$
and let $N^*_{\eta_i}\preceq M$ be any $\ae$-prime model over
$N^*_{\eta_i^-}\cup N_{\eta_i}$.
Then it is easily checked that  $\d^*=\<N^*_\eta,a_\eta:\eta\in I\>$
is an $\aeP$-decomposition inside $M$ that is a blow up of $\d$.
\endproof

`Blowing down' a decomposition is more delicate and requires  two 
technical Lemmas, 
Lemma~\ref{oneone} and Lemma~\ref{twotwo} that 
assert the existence of $\ae$-submodels of a given $\ae$-saturated structure
with certain properties.

\begin{Lemma}  \label{blowdown}  Suppose that $M$ is $\ae$-saturated
and $\d=\<N_\eta,a_\eta:\eta\in I\>$ is a prime $\aeP$-decomposition
inside $M$.  For any  $\ae$-saturated $N_0\preceq N_{\<\>}$ such that for
every $\eta\in Succ_I(\<\>)$, either $\tp(a_\eta/N_{\<\>})$ does not
fork over $N_0$ or $\tp(a_\eta/N_{\<\>}$ is regular and non-orthogonal to $N_0$.  Then there is a prime $\aeP$-decomposition
$\d_0$ inside $M$ with $N^{\d_0}_{\<\>}=N_0$  such that 
$\d$ is a blow up of $\d_0$.
\end{Lemma}

\bp  Choose an enumeration $\<\eta_i:i<i^*\>$ of $I$ as in the proof
of Lemma~\ref{blowup}.  That is, $\eta_0=\<\>$, $\eta_i\triangleleft\eta_j$
implies $i<j$, and $\eta_i\in Succ_I(\<\>)$ if and only if
$1\le i<\alpha^*$ for some $\alpha^*\le i^*$.

Put $N^0_{\eta_0}=N_0$.  For $1\le i<i^*$ 
we inductively construct $N^0_{\eta_i}$ to satisfy:
\begin{itemize}
\item  $N^0_{\eta_i}\preceq N_{\eta_i}$ and 
$\fg {N^0_{\eta_i}} {N^0_{(\eta_i^-)}} {N^1_{\eta_i^-}}$
\item  $N_{\eta_i}$ is $\ae$-prime over $N_{\eta_i}^0\cup N_{\eta_i^-}$
\item  $a_{\eta_i}\in N^0_{\eta_i}$ and $N^0_{\eta_i}$ is $\ae$-prime over
$N^0_{\eta_i}\cup\{a_{\eta_i}\}$.
\end{itemize}

To accomplish this,
for each $1\le i<\alpha^*$, use Lemma~\ref{twotwo}
to define $N_{\eta_i}^0$ (where $M_1=N_0$, $M=N_{\eta_i}$).  We can
take $N_{\eta_i}^0$ to be the $N$ there, and we can take $a^*$ to be
$a_{\eta_i}$.  Similarly, for $\alpha^*\le i<i^*$ we apply Lemma~\ref{oneone},
where $M$ is taken to be $N_{\eta_i^-}$, $M_1$ is $N^0_{\eta_i^-}$, $M_2$ is $N_{\eta_i^{--}}$,
$a$ is $a_{\eta_i}$,
and taking $N^0_{\eta_i}$ to be the $N$ produced there.
\endproof

\begin{Definition}  {\em Suppose $M$ is $\ae$-saturated and 
${B\choose A}\in\Gamma_\P(M)$.  We say that an $\eps$-finite subset $W\subseteq M$
{\em has a base $W_0\subseteq W$ respecting  $B\choose A$} if
 $A\subseteq W_0$, $\fg {W_0} A B$,
and $W$ is dominated by $B$ over $W_0$.
}
\end{Definition}

\begin{Lemma}  \label{getpartition}
If $\d=\<N_\eta,a_\eta:\eta\in I\>$ is a $\aeP$-decomposition inside $M$ over $B\choose A$
and $V\subseteq \bigcup_{\eta\in I} N_\eta$ is $\eps$-finite, then
there is an $\eps$-finite $W$ with $V\subseteq W$ and
$W\setminus V\subseteq N_{\<\>}$
that has a base $W_0\subseteq W\cap N_{\<\>}$ respecting $B\choose A$.
\end{Lemma}

\bp    Without loss, we may assume $A\subseteq V$.
It follows from the definition of an $\aeP$-decomposition inside $M$
over $B\choose A$ that $\fg B A {N_{\<\>}}$
and that  $B$ dominates $\bigcup_{\eta\in I} N_\eta$ 
and hence $V$ over $N_{\<\>}$. As both $B$ and $V$ are $\eps$-finite,
it follows from superstability that there is an $\eps$-finite $C\subseteq N_{\<\>}$
such that $\fg {BV} C {N_{\<\>}}$.  So $B$ dominates $V$ over $C$.  
Again, without loss, $A\subseteq C$.
Take $W=V\cup C$.  Then $W_0:=W\cap N_{\<\>}$ is a base respecting
$B\choose A$.
\endproof

\begin{Lemma} \label{twopiece}
 Suppose $W\subseteq M$ is $\eps$-finite and has
a base $W_0\subseteq W$ respecting $B\choose A$ and that $N\preceq M$ is
$\ae$-prime over $\emptyset$, $W_0\subseteq N$, with $\fg N A B$.
Then
there is $N[B]\preceq M$ that is $\ae$-prime over $N\cup B$
such that
$W\subseteq N[B]$ and such that the two-element sequence
$\<N,N[B]\>$ is an $\aeP$-decomposition inside $M$ over $B\choose A$
(taking $a_{\<0\>}$ to be $B$).
\end{Lemma}

\bp  As $A,B,W$ are all $\eps$-finite, $\fg N A B$, $W$ dominated
by $B$ over $W_0$,
and the fact that $\tp(W/\acl(W_0\cup B))$ is stationary,
 it follows that $\tp(W/NB)$ is $\aleph_\eps$-isolated.
Thus, $W\subseteq N[B]$ for some $\ae$-prime model over  $N\cup B$.
Checking that the two element sequence is an $\aeP$-decomposition inside
$M$ over $B\choose A$ is routine.
\endproof

\begin{Proposition}  \label{same}
If $M$ is $\ae$-saturated, ${B\choose A}\in\Gamma_\P(M)$, and the
$\eps$-finite set $W\subseteq M$ has a base $W_0\subseteq W$ respecting 
$B\choose A$,
then there is an $\aeP$-decomposition $\d$ of $M$ over $B\choose A$
with  $W_0\subseteq N^\d_{\<\>}$ and $W\subseteq N^\d_{\<0\>}$.
Moreover, if $\d_0$ is any $\aeP$-decomposition of $M$ over $B\choose A$,
then $\d$ can be chosen so that $\PP(\d,M)=\PP(\d_0,M)$.
\end{Proposition}

\bp  Suppose $\d_0=\<N_\eta,a_\eta:\eta\in I\>$ is given.  As 
$\dcl(a_{\<0\>})=\dcl(B)$, we may assume that $a_{\<0\>}=B$.
Thus, $\fg B A {N_{\<\>}}$.
Choose a finite $D\subseteq N_{\<\>}$ such that $A\subseteq D$
 and $\fg W {DB} {N_{\<\>}B}$.
By e.g., 1.18(9) of \cite{Sh401} there is $N^1\preceq N_{\<\>}$
that is $\ae$-prime over $\emptyset$, $\fg {N^1} A D$, and
$N_{\<\>}$ is $\ae$-prime over $N^1\cup D$.

As for non-forking, we claim that the set $\{B,W_0,N^1\}$ is independent
over $A$.  To see this, first recall that $\fg B A {W_0}$ since $W_0$ is 
a base respecting $B\choose A$.  As well,
$\fg B A {N_{\<\>}}$ since $\d$ is over $B\choose A$.  Thus,
$\fg B D {N_{\<\>}}$.  Thus, by our choice of $D$ and forking calculus,
$\fg {W_0B} D {N_{\<\>}}$, so $\fg {W_0B} D {N^1}$ since 
$N^1\preceq N_{\<\>}$.  But now, as $\fg {N^1} A D$, we have
$\fg {N^1} A {BW_0}$ which gives the independence.

By Lemma~\ref{blowdown}, there is an $\aeP$-decomposition
$\d_1$ inside $M$ over $B\choose A$ with $N^{\d_1}_{\<\>}=N^1$.
By Lemma~\ref{1.20}(4) $\d_1$ is an $\aeP$-decomposition of $M$ over
$B\choose A$, so by Lemma~\ref{blowupPP} $\PP(\d_0,M)=\PP(\d_1,M)$.

Next, let $N^2\preceq M$ be $\ae$-prime over $N^1\cup W_0$.
As $\fg B {N^1} {W_0}$, and the $\ae$-isolation of $N^1$ we have 
$\fg {N^2} {N^1} B$.  Thus, by Lemma~\ref{blowup} there is
a $\aeP$-decomposition $\d_2$ inside $M$ over ${B\choose A}$ with 
$N^{\d_2}_{\<\>}=N^2$.  Again, by Lemma~\ref{1.20}(4)
$\d_2$ is an $\aeP$-decomposition of $M$ over $B\choose A$
and by Lemma~\ref{blowupPP} $\PP(\d_1,M)=\PP(\d_2,M)$.

Put $N:=N^2$.
Clearly, $W_0\subseteq N$ and we showed $\fg B {N^1} {N}$.
But, as $B$ and $N^1$ are independent over $A$, $\fg B A {N}$.
So, by Lemma~\ref{twopiece} there is $N[B]\preceq M$, 
$\ae$-prime over $N\cup B$, such that $W\subseteq N[B]$ and
$\<N,N[B]\>$ is an
$\aeP$-decomposition inside $M$ over $B\choose A$.

Finally, by Lemma~\ref{straightBA} there is an $\aeP$-decomposition
$\d_3$ of $M$ over $B\choose A$ end extending $\<N,N[B]\>$.
As $N^{\d_3}_{\<\>}=N=N^{\d_2}_{\<\>}$ we conclude by Lemma~\ref{charPP}
that $\PP(\d_3,M)=\PP(\d_2,M)$.  Thus, $\PP(\d_3,M)=\PP(\d_0,M)$
and we finish.
\endproof

We are finally ready to prove our main Theorem.

\begin{Theorem}  \label{big}
Suppose that $M$ is $\ae$-saturated and ${B\choose A}\in\Gamma_\P(M)$.
Then $\PP(\d_1,M)=\PP(\d_2,M)$ for any two prime $\aeP$-decompositions
$\d_1,\d_2$ of $M$ over $B\choose A$.
\end{Theorem}

\bp  Suppose $\d_1=\<N_\eta,a_\eta:\eta\in I\>$.  By symmetry, it suffices
to prove that every $p\in\PP(\d_1,M)$ is in $\PP(\d_2,M)$.
Fix such a $p$ and choose $\eta\in I\setminus\{\<\>\}$ such that
$p\not\perp N_\eta$ but $p\perp N_{\eta^-}$.  Choose $q\in S(N_\eta)$
regular such that $p\not\perp q$ and choose a finite $V\subseteq N_\eta$
on which $q$ is based and stationary.  By Lemma~\ref{getpartition}
there is an $\eps$-finite $W$ such that $V\subseteq W\subseteq M$
that has a subset $W_0=W\cap N_{\<\>}$ respecting $B\choose A$.
Note that since $p\perp N_{\<\>}$ we have $p\perp W_0$, hence $q\perp W_0$.
 By applying Proposition~\ref{same} to $W$ and
$\d_2$, we get that there is a prime $\aeP$-decomposition $\d^*$ of $M$ 
over $B\choose A$ with $\PP(\d^*,M)=\PP(\d_2,M)$.  But, by construction,
there is a type parallel to $q$ (and hence non-orthogonal to $p$) in $S(N^{\d_2}_{\<0\>})$. 
As well, since $B$ dominates $W$ over $W_0$ and $\fg B A {N_{\<\>}}$
we have $\fg W {W_0} {N_{\<\>}}$.
As $q$ is based on $W$ and $q\perp W_0$, we have that $q$ (and hence $p$)
is orthogonal to $N_{\<\>}$. Thus, $p\in\PP(\d^*,M)=\PP(\d_2,M)$.
\endproof

The previous Theorem inspires the following definition.

\begin{Definition}  {\em
 For $M$ $\ae$-saturated and ${B\choose A}\in\Gamma_\P(M)$,
$\PP({B\choose A},M)=\PP(\d,M)$ for some (equivalently for every)
prime $\aeP$-decomposition $\d$ of $M$ over $B\choose A$.
}
\end{Definition}

\begin{Corollary}  \label{fulltree}
Suppose that $M$ is $\ae$-saturated, ${B\choose A}\in\Gamma_\P(M)$,
and that $\d=\<N_\eta,a_\eta:\eta\in I\>$ is a prime, $\aeP$-decomposition
of $M$ satisfying (1) $N_{\<\>}$ is $\ae$-prime over $A$; (2) $\fg B A {N_{\<\>}}$; and (3) $N_{\<0\>}$ is $\ae$-prime over $N_{\<\>}\cup B$.
Then, for every $p\in S(M)\cap\P$, $p\in \PP({B\choose A},M)$ if and only
if $\<0\>\trianglelefteq \eta(p)$, where $\eta(p)$ is the unique 
$\trianglelefteq$-minimal $\eta\in I$ satisfying $p\not\perp N_\eta$ (see
Corollary~\ref{mineta}).
\end{Corollary}

\bp  Given $\d=\<N_\eta,a_\eta:\eta\in I\>$ as above, let
$X=\{\nu\in I\setminus\{\<\>\}:\neg(\<0\>\trianglelefteq\nu\}$
and let $I_0=I\setminus X$.  The conditions on $\d$ ensure that
$\d_0:=\<N_\eta,a_\eta:\eta\in I_0\>$ is a prime, $\aeP$-decomposition
of $M$ above $B\choose A$.  Thus, by Theorem~\ref{big}, for any
$p\in S(M)\cap \P$ we have
$$p\in\PP({B\choose A},M)\quad\Leftrightarrow\quad p\in\PP(\d_0,M)\quad
\Leftrightarrow\quad \<0\>\trianglelefteq \eta(p)$$
\endproof

The following characterization is analogous to Claim~1.24  of \cite{Sh401}.

\begin{Proposition}  Assume that $M_1\preceq M_2$ are $\ae$-saturated
and ${B\choose A}\in\Gamma_\P(M_1)$.  Then the following are equivalent:
\begin{enumerate}
\item  No $p\in\PP({B\choose A},M_1)$ is realized in $M_2$;
\item  There is a prime $\aeP$-decomposition of $M_1$ above $B\choose A$
that is also a prime $\aeP$-decomposition of $M_2$ above $B\choose A$; and
\item  Every  prime $\aeP$-decomposition of $M_1$ above $B\choose A$
 is also a prime $\aeP$-decomposition of $M_2$ above $B\choose A$.
\end{enumerate}
\end{Proposition}

\bp  $(3)\Rightarrow(2)$ is immediate since prime $\aeP$-decompositions
of $M_1$ over $B\choose A$ exist.

$(2)\Rightarrow(1)$:  Let $\d=\<N_\eta,a_\eta:\eta\in I\>$ be a 
prime $\aeP$-decomposition of 
$M_1$ above $B\choose A$ 
and  assume that there
is $e\in M_2\setminus M_1$ such that $p=\tp(e/M_1)\in\PP(\d,M_1)$.
Choose $\eta\in I$ to be $\triangleleft$-minimal such that $p\not\perp N_\eta$.
Note that $\<0\>\trianglelefteq \eta$.
By Fact~\ref{refer}(4) and because $N_\eta,M_1,M_2$ are
$\aleph_\eps$-saturated, we can replace $e$ by the realization
of a non-orthogonal regular type that satisfies $\fg e {N_\eta} {M_1}$.
As $e\in C_\eta(M_2)$, $\{a_\nu:\nu\in Succ_I(\eta)\}$ is not maximal
$N_\eta$-independent subset of $C_\eta(M_2)$, so $\d$ is not a prime $\aeP$-decomposition of 
$M_2$ above $B\choose A$.

$(1)\Rightarrow(3)$: 
Let $\d=\<N_\eta,a_\eta:\eta\in I\>$ be a prime $\aeP$-decomposition of 
$M_1$ above $B\choose A$, and assume that it is not a prime $\aeP$-decomposition
of $M_2$.  Then, by Definition~\ref{ofMoverBA}, there is $\eta\in I\setminus\{\<\>\}$
such that $\{a_\nu:\nu\in Succ_I(\eta)\}$ is not a maximal, $N_\eta$-independent subset
of $C_\eta(M_2)$.  As both $N_\eta$ and $M_2$ are $\ae$-saturated, this implies that there is $e\in C_\eta(M_2)$
such that $\tp(e/N_\eta)\in\P$, but $\fg e {N_\eta} {\{a_\nu:\nu\in Succ_I(\eta)\}}$.
By Lemma~\ref{straightBA}(5), $\d$ $\P$-exhausts $M_1$ over $B\choose A$ so
$\fg e {N_\eta} {M_1}$.
Thus, $p=\tp(e/M_1)$ is an element of 
$\PP({B\choose A},M_1)$ that is realized in $M_2$.
\endproof

For pairs ${B_1}\choose {A_1}$ and ${B_2}\choose {A_2}$ from $\Gamma(M)$,
we consider two ways in which ${B_2}\choose {A_2}$ can extend 
${B_1}\choose{A_1}$, corresponding to the former having `more information'
or `appearing higher up in a $\P$-decomposition.'

First, write ${{B_1}\choose {A_1}}\le_a{{B_2}\choose {A_2}}$ if both are from
$\Gamma(M)$, $A_1\subseteq A_2$, $B_1\subseteq B_2$, 
$\fg {B_1} {A_1} {A_2}$, and $B_2$ dominated by $B_1$ over $A_2$.
Intuitively, think of ${B_2}\choose {A_2}$ as being a `better approximation'
of $(N,N',a)$.  

The next approximation, which should be thought of as `stepping up in the tree'
is given by 
${{B_1}\choose {A_1}}\le_b{{B_2}\choose {A_2}}$ if and only if
$A_2=B_1$, and $\tp(B_2/A_2)$ is regular and is orthogonal to $A_1$.

Finally, let $\le^*$ be the transitive closure of $\le_a\cup\le_b$.

\begin{Proposition}  \label{compare}
Fix an $\ae$-saturated model $M$ and  ${B_1}\choose {A_1}$, 
${B_2}\choose {A_2}$ from $\Gamma_\P(M)$.
\begin{enumerate}
\item  If ${{B_1}\choose{A_1}}\le_a {{B_2}\choose{A_2}}$, then
$\PP({{B_1}\choose{A_1}},M)=\PP({{B_2}\choose{A_2}},M)$;

\item  If ${{B_1}\choose{A_1}}\le_b {{B_2}\choose{A_2}}$, then
$\PP({{B_2}\choose{A_2}},M)$ is a proper subset of
$\PP({{B_1}\choose{A_1}},M)$ 

\item  If ${{B_1}\choose{A_1}}\le^* {{B_2}\choose{A_2}}$ then
$\PP({{B_2}\choose{A_2}},M)\subseteq
\PP({{B_1}\choose{A_1}},M)$;

\item  If $A_1=A_2$ (whose common value we denote by $A$)
  $\tp (B_1/A)$,  $\tp(B_2,A)$ are both
regular, and $\nfg {B_1} A {B_2}$, then
$\PP({{B_1}\choose{A}},M)=\PP({{B_2}\choose{A}},M)$.

\item If $A_1=A_2=A$ and $\fg {B_1} A {B_2}$, then the sets
$\PP({{B_1}\choose A},M)$ and $\PP({{B_2}\choose A}, M)$
are disjoint.
\end{enumerate}
\end{Proposition}

\bp  (1)  Let $N_{\<\>}\preceq M$ be $\aleph_\eps$-prime over $\emptyset$ with
$A_2\subseteq N_{\<\>}$ and $\fg {B_2} {A_2} {N_{\<\>}}$.  Let $N_{\<\>}$ be
$\ae$-prime over $N_{\<\>}\cup B_2$, let $a_{\<0\>}=B_2$,
and let $\d=\<N_\eta,a_\eta:\eta\in I\>$ be a prime $\aeP$-decomposition of $M$ over
${{B_2}\choose {A_2}}$ end extending $\<N_{\<\>},N_{\<0\>}\>$.
It follows easily by the forking calculus that $\d$ is also a prime $\aeP$-decomposition
of $M$ over
${{B_1}\choose {A_1}}$.  Thus, two applications of Theorem~\ref{big} yield
$$\PP({{B_2}\choose {A_2}},M)=\PP(\d,M)=\PP({{B_1}\choose{A_1}},M)$$
  
(2)  Given $A_1\subseteq B_1=A_2\subseteq B_2\subseteq M$ with $\tp(B_2/A_2)\perp A_1$,
first choose an $\ae$-prime $N_{\<\>}\preceq M$ containing $A_1$ with $\fg {B_2} {A_1} {N_{\<\>}}$.
Note that $\tp(B_2/A_2)\perp N_{\<\>}$.  Let $a_{\<\>}$ be an arbitrary element of $N_{\<\>}$, let
$a_{\<0\>}:=A_2$, and choose $N_{\<0\>}\preceq M$ to be $\ae$-prime over $N_{\<\>}\cup A_2$, with
$\fg {N_{\<0\>}} {N_{\<\>}A_2} {B_2}$.  Also, choose $N_{\<0,0\>}\preceq M$ to be $\ae$-prime over
$N_{\<0\>}\cup B_2$ and let $a_{\<0,0\>}:=B_2$.

Let $\d_0$ be the three-element prime $\aeP$-decomposition $\<N_\eta,a_\eta:\eta\in\{\<\>,\<0\>,\<0,0\>\}\>$
inside $M$ above ${{B_1}\choose {A_1}}$.  
Next, by `collapsing', let $\d_0'=\<N_\eta',a_\eta':\eta\in\{\<\>,\<0\>\}\>$ be the two-element
prime $\aeP$-decomposition inside $M$
above ${{B_2}\choose {A_2}}$, where $N'_{\<\>}:=N_{\<0\>}$, $a'_{\<\>}:=a_{\<0\>}$, 
$N'_{\<0\>}:=N_{\<0,0\>}$, and $a'_{\<0\>}:=a_{\<0,0\>}$.  

Next, choose a prime $\aeP$-decomposition $\d'=\<N_\eta',a_\eta':\eta\in I'\>$ of $M$ above ${{B_2}\choose {A_2}}$
end extending $\d_0'$.  It follows immediately from Theorem~\ref{big} that 
$\PP({{B_2}\choose {A_2}},M)=\PP(\d',M)$, so to obtain the inclusion
$\PP({{B_2}\choose{A_2}},M)\subseteq
\PP({{B_1}\choose{A_1}},M)$ it suffices to construct a prime $\aeP$-decomposition $\d=\<N_\eta,a_\eta:\eta\in J\>$ 
inside $M$ over
${{B_1}\choose {A_1}}$
such that, for any $p\in S(M)\cap \P$, if $p\not\perp N_\eta'$ but $p\perp N_{\eta^-}'$ for some $\eta\in I'$
with $\<0\>\trianglelefteq \eta'$, there is $\eta\in J$ such that $\<0\>\trianglelefteq\eta$, $p\not\perp N_\eta$,
but $p\perp N_{\eta^-}$.  

We accomplish this as follows: Recall that $N_\eta,a_\eta$ were defined for $\eta\in\{\<\>,\<0\>,\<0,0\>\}$
above.  Let $J'\subseteq I'$ be $\{\<\>\}\cup\{\eta\in I':\<0\>\trianglelefteq\eta\}$, and define a function
$h$ with domain  $J'$ by
$h(\eta):=\<0\>\conc\eta$ if $\eta\neq\<\>$.
That is, the function $h$ is `undoing' the collapse given above.
Let $J=\{\<\>,\<0\>\}\cup\{h(\eta:\eta\in J'\})$, and for each $\eta\in J'$,
put $N_{h(\eta)}:=N'_{\eta}$ and $a_{h(\eta)}:=a_\eta'$.
Then $\d:=\<N_\eta,a_\eta:\eta\in J\>$ is a prime $\aeP$-decomposition inside $M$ above 
${{B_1}\choose {A_1}}$, and for any $p\in \PP(\d',M)$, if $p\not\perp N'_\eta$ for some
$\eta\in J'$, then $p\not\perp N_{h(\eta)}$.  Thus, $\d$ is as required.

To show that the inclusion is strict, choose any regular type $q\in S(N_{\<0\>})$ that is non-orthogonal to
$\tp(B_2/N_{\<0\>})$.  It is easy to check that the non-forking extension of $q$ to $S(M)$ is an element of
$\PP({{B_1}\choose{A_1}},M)\setminus
\PP({{B_2}\choose{A_2}},M)$.

(3) follows immediately from (1) and (2).

(4)  By symmetry, it suffices to show that
$\PP({{B_2}\choose A},M)\subseteq\PP({{B_1}\choose A},M)$, so fix a 
regular type $p\in S(M)\cap\P\setminus\PP({{B_1}\choose A},M)$.
We will eventually produce a prime $\aeP$-decomposition $\d_2$ inside $M$ over
${{B_2}\choose A}$ with the property that $p\not\perp N_\eta$ for some
$\eta$ satisfying $\neg(\<0\>\trianglelefteq\eta)$, which suffices by Lemma~\ref{straightBA}(4) and Theorem~\ref{big}.

We begin by choosing an $\ae$-prime (over $\emptyset$) $N_{\<\>}\preceq M$
that contains $A$, but $\fg {B_1B_2} A {N_{\<\>}}$.  Note that
$B_1$ and $B_2$ are domination equivalent over $N_{\<\>}$.

Let $a_{\<\>}\in N_{\<\>}$
be arbitrary, let $N^1$ be $\ae$-prime over $N_{\<\>}\cup B_1$, and let $a_{\<0\>}:=B_1$.
Then $\d_1:=\<N_\eta,a_\eta:\eta\in\{\<\>,\<0\>\}\>$ is a two-element prime $\aeP$-decomposition
inside $M$ over ${{B_1}\choose A}$.  Let $\d_1'=\<N_\eta,a_\eta:\eta\in I\>$
be a prime $\aeP$-decomposition of $M$ over ${{B_1}\choose A}$ end extending $\d_1$.
Next, let $\d_1^*=\<N_\eta,a_\eta:\eta\in J\>$ be a prime $\aeP$-decomposition of $M$
end extending $\d_1'$.  
Let $H=\{\eta\in J:\neg(\<0\>\trianglelefteq\eta)\}$.  Then
$H$ is a subtree of $J$, whose intersection with $I$ is $\{\<\>\}$.
Furthermore, as $p\in\P$, it follows from Corollary~\ref{mineta} that 
$p\not\perp N_\eta$ for some $\eta\in J$.  However, since $p\not\in\PP({{B_1}\choose A},M)$,
it follows from Theorem~\ref{big} that $p\not\in\PP(\d_1',M)$, hence
$p\not\perp N_\eta$ for some $\eta\in H$.

But now, choose $N^2\preceq M$ to be $\ae$-prime over $N_{\<\>}\cup B_2$.
Let $\d_2:=\<N_\eta,a_\eta:\eta\in H\>\conc(N^2,B_2)$.  As $B_1$ and $B_2$
are domination equivalent over $N_{\<\>}$, it is easily checked that
$\d_2$ is a prime $\aeP$-decomposition inside $M$ over ${{B_2}\choose A}$.
Let $\d_2^*=\<N_\eta,a_\eta:\eta\in I_2\>$ 
be any prime $\aeP$-decomposition inside $M$ end extending $\d_2$.
But, as $p\not\perp N_\eta$ for some $\eta\in H$, it follows from independence
that $p\perp N_\nu$ for any $\nu\in I_2$ satisfying $\<0\>\trianglelefteq \nu$.
Thus, $p\not\in \PP({{B_2}\choose A},M)$ by Theorem~\ref{big} again.

(5)  Let $N_{\<\>}\preceq M$ be $\ae$-prime over $A$ with 
$\fg {N_{\<\>}} A {B_1B_2}$ and choose an $\eps$-finite $B_0\in N_{\<\>}$
arbitrarily.  For $\ell=1,2$, choose $N_{\<\ell\>}$ to be $\ae$-prime
over $N_{\<\>}\cup B_\ell$.  Clearly, 
$$\d':=\{(N_{\<\>},B_0),(N_{\<0\>},B_1),(N_{\<1\>},B_1)\}$$ is a three
element, prime $\aeP$-decomposition
inside $M$.   By Lemma~\ref{maxendextension}(2) there is a prime, 
$\aeP$-decomposition
$\d=\<N_\eta,a_\eta:\eta\in I\>$ of $M$ end extending $\d'$.
It is easily checked that $\d$ satisfies the hypotheses of Corollary~\ref{fulltree},
as does the modification formed by exchanging the roles of $\<0\>$ and $\<1\}$.
Thus, for any $p\in S(M)\cap \P$, we have
$p\in\PP({{B_1}\choose A},M)$ if and only if $\<0\>\trianglelefteq \eta(p)$,
and that $p\in\PP({{B_2}\choose A}, M)$ if and only if $\<1\>\trianglelefteq
\eta(p)$.  As the elements $\<0\>$ and $\<1\>$ are incompatible, it
follows that the sets $\PP({{B_1}\choose A},M)$ and $\PP({{B_2}\choose A}, M)$
are disjoint.
\endproof

\begin{Corollary} \label{incomparable}
 Suppose that $M$ is $\aleph_\epsilon$-saturated and that
$\d=\<N_\eta,a_\eta:\eta\in I\>$ is any weak $\P$-decomposition inside $M$.
Choose any incomparable nodes $\eta_1,\eta_2\in I$.  If, for each $\ell=1,2$,
$A_\ell\subseteq N_{\eta_\ell}^-$ is $\epsilon$-finite on which 
$\tp(a_{\eta_\ell}/N_{\eta_\ell}^-)$ is based and stationary and
$B_\ell=\acl(A_\ell\cup\{a_{\eta_\ell}\})$, then
the sets $\PP({{B_1}\choose{A_1}},M)$ and
$\PP({{B_2}\choose{A_2}},M)$ are disjoint.
\end{Corollary}

\bp  As $\eta_1$ and $\eta_2$ are incomparable, neither one is $\<\>$,
so let $\mu$ denote the meet $\eta_1^-\wedge\eta_2^-$.  By incomparability
again, there are distinct ordinals $\alpha_1\neq\alpha_2$ such that
$\mu\conc\<\alpha_1\>\trianglelefteq\eta_1$, while 
$\mu\conc\<\alpha_2\>\trianglelefteq\eta_2$.
Choose an $\epsilon$-finite $E\subseteq M_{\mu}$ over which both
types $\tp(a_{\mu\conc\<\alpha_\ell\>}/M_\mu)$ are based and stationary,
and let $C_\ell=\acl(a_{\mu\conc\<\alpha_\ell\>}\cup E)$ for each $\ell$.
As $\fg {C_1} E {C_2}$ it follows from Proposition~\ref{compare}(5) that 
the sets $\PP({{C_1}\choose{E}},M)$ and
$\PP({{C_2}\choose{E}},M)$ are disjoint.  But, by Proposition~\ref{compare}(3)
$\PP({{B_\ell}\choose{A_\ell}},M)\subseteq
\PP({{C_1}\choose{E}},M)$ for each $\ell$ and the result follows.
\endproof

\begin{Proposition}  \label{trivialapp}
Suppose that $M$ is $\ae$-saturated and $p_1\in S(A_1)$, $p_2\in S(A_2)$
are non-orthogonal, trivial, regular types over $\eps$-finite subsets of $M$.
If, for $\ell=1,2$, $I_\ell$ is a maximal, $A_\ell$-independent subset of
$p_\ell(M)$, then there are cofinite subsets $J_\ell\subseteq I_\ell$ and
a bijection $h:J_1\rightarrow J_2$ such that
$$\PP({c\choose {A_1}},M)=\PP({{h(c)}\choose {A_2}},M)$$
for every $c\in J_1$.
\end{Proposition}

\bp  Let $D=A_1\cup A_2$.  For $\ell=1,2$, let
$J_\ell:=\{c\in I_1:\fg c {A_\ell} D\}$ and let $q_\ell$ denote the non-forking
extension of $p_\ell$ to $S(D)$.  
Then $J_\ell$ is a cofinite subset of $I_\ell$ and is a maximal, $D$-independent
subset of $q_\ell(M)$.  As the regular types are trivial and non-orthogonal,
$p_1$ and $p_2$ are not almost orthogonal, so as $M$ is $\ae$-saturated,
we have that for every $c\in q_1(M)$, there is $c'\in q_2(M)$ such that
$\nfg {c_1} D {c_2}$.   It follows that there is a unique bijection
$h:J_1\rightarrow J_2$ satisfying $\nfg c D {h(c)}$ for each $c\in J_1$.
Thus, $\PP({c\choose {A_1}},M)=\PP({{h(c)}\choose {A_2}},M)$
by Clauses (1) and (4) of  Proposition~\ref{compare}.
\endproof

\section{Decompositions and non-saturated models}

Until this point, we have been looking at various flavors of decompositions of
$\aleph_\eps$-saturated models.  It would be desirable to see what effect
these results have on understanding decompositions of arbitrary models.
In the first subsection, given an arbitrary model $M$ and a sufficiently
saturated elementary extension $M^*$, one can produce an $\aeP$-decomposition
$\d=\<M_\eta,a_\eta:\eta\in I\>$
of $M^*$ that `enumerates $M$ as slowly as possible.'  In particular, given any
$\eps$-finite $A\subseteq M$, there is a finite subtree $J\subseteq I$, an elementary
submodel $M^J\preceq M^*$ that is $\aleph_\eps$-prime over $\bigcup_{\eta\in J} M_\eta$,
and an $\eps$-finite $B$, $A\subseteq B\subseteq M$ that satisfy
$\fg M B {M_J}$.  

In the second subsection, we obtain a weak uniqueness result for $\P$-decompositions of
unsaturated models $M$ satisfying certain constraints.  Whereas these conditions seem contrived,
Theorem~\ref{quotable} plays a major role in \cite{LSh2}.

\subsection{Large extensions of weak decompositions}

As usual, we assume that $\P$ is a set of stationary, regular types closed under isomorphism and  non-orthogonality,
and we assume that our theory $T$ is superstable with $\P$-NDOP.

\begin{Definition}  \label{respectingM}  
{\em  Suppose that $M\preceq M^*$ are given, with $M$ arbitrary, but $M^*$ sufficiently saturated.
A prime $\aeP$-decomposition $\d^*=\<N_\eta,a_\eta:\eta\in I\>$ of $M^*$ {\em respects $M$}
if
there is a continuous, elementary chain $\<M_\alpha:\alpha\le\alpha^*\>$ of $\aleph_\eps$-saturated
elementary substructures of $M^*$ with $\bigcup_{\alpha\le\alpha^*} M_\alpha=M^*$;
a sequence $\<\d_\alpha:\alpha\le\alpha^*\>$ of prime $\aeP$-decompositions of $M_\alpha$ 
with $\d_{\alpha^*}=\d^*$; and
a sequence $\<a_\alpha:\alpha\le\alpha^*\>$ of elements from $M^*$ that satisfy the following constraints:
\begin{enumerate}
\item  $M_0=N_{\<\>}$ and the  sets $M$ and $M_0$ are independent;
\item  If $\beta\le\alpha$ then $\d_\alpha$ end extends $\d_\beta$ with $\d_\gamma=\bigcup\d_\alpha$ for $\gamma$ a limit ordinal;
\item  The trees $I_\alpha$ indexing the decompositions $\d_\alpha$ satisfy $|I_{\alpha+1}\setminus I_\alpha|\le 1$ for each $\alpha<\alpha^*$;
\item  If $I_{\alpha+1}\setminus I_\alpha=\{\eta\}$, then $N_\eta$ is $\aleph_\eps$-prime over $N_{\eta^-}\cup\{a_\alpha\}$
and $\fg {N_\eta} {N_{\eta^-}a_\alpha} {MM_\alpha}$.
\end{enumerate}
}
\end{Definition}

\begin{Lemma}  \label{respectingMexistence}
Suppose that $M\preceq M^*$, where $M^*$ is saturated and $||M^*||>||M||+2^{|T|}$.
Then a prime $\aeP$-decomposition $\d^*$ of $M^*$ respecting $M$ exists.
\end{Lemma}

\bp  We recursively construct sequences $\<M_\alpha\>$, $\<\d_\alpha\>$ and $\<a_\alpha\>$ with
the additional constraint of $||M_\alpha||<||M^*||$ for each $\alpha<\alpha^*$ as follows.
First, choose $N_{\<\>}\preceq M^*$ to be $\aleph_\eps$ saturated with $\fg {N_{\<\>}} {\phantom{X}} M$,
let $M_0=N_{\<\>}$, $I_0=\{\<\>\}$, and $\d_0=\<N_{\<\>}\>$.
For $\alpha\le\alpha^*$ a limit ordinal, simply take unions.

Next, fix an enumeration $\<c_i:i<\lambda\>$ of $M^*$ with $\lambda=||M^*||$ and
the elements of $M$ forming an initial segment
and assume that $M_\beta$ and  $\d_\beta$ have been defined.  Let $c^*$ be the least element of $M^*$ that is
not an element of $M_\beta$.  There are now two cases, depending on $\tp(c^*/M_\beta)$.

\medskip

\par\noindent{\bf Case 1:}  $\tp(c^*/M_\beta)\perp \Pactive$.

\medskip  In this case, choose a regular type $q\in S(M_\beta)$ non-orthogonal to $\tp(c^*/M_\beta)$.
As $M^*$ is saturated, choose an element $a_\beta\in M^*$ realizing $q$ with $\nfg {a_\beta} {M_\beta} {c^*}$.
Let $\d_{\beta+1}=\d_{\beta}$, and let $M_{\beta+1}\preceq M^*$ be $\aleph_\eps$-prime over 
$M_\beta\cup\{a_\beta\}$ and satisfying $\fg {M_{\beta+1}} {M_\beta a_\beta}
M$.

\medskip

\par\noindent{\bf Case 2:}  $\tp(c^*/M_\beta)\not\perp \Pactive$.

\medskip  In this case, choose a regular type $q\in S(M_\beta)\cap\Pactive$ non-orthogonal to $\tp(c^*/M_\beta)$.
By Corollary~\ref{mineta}, there is a unique $\eta\in I_\beta$ such that $q\not\perp N_\eta$, but $q\perp N_{\eta^-}$
(if $\eta\neq\<\>$).  Without loss, we may assume that $q$ does not fork over $N_\eta$.  As $M^*$ is saturated,
we can choose an element
$a_\beta\in M^*$ realizing $q$ with $\nfg {a_\beta} {M_\beta} {c^*}$.
Let $\gamma$ be the least ordinal such that $\nu:=\eta\conc\<\gamma\>\not\in I_\beta$.
Choose $N_\nu\preceq M^*$ to be $\aleph_\eps$-prime over $N_\eta\cup\{a_\beta\}$ and satisfying $\fg {N_\nu} {N_\eta\cup\{a_\beta\}}
{MM_\beta}$.  As $N_\eta$
is $\aleph_\eps$-saturated, it follows by Fact~\ref{refer}(2)
that $\fg {N_\nu} {N_\eta} {M_\beta}$.  Choose $M_{\beta+1}\preceq M^*$ to be $\aleph_\eps$-prime
over $M_\beta\cup N_\nu$ and satisfying $\fg {M_{\beta+1}} {M_\beta N_\nu}
M$.
Let $I_{\beta+1}=I_\beta\cup\{\nu\}$, let 
Let $\d_{\beta+1}=\d_\beta\conc\<N_\nu\>$, and let $M_{\beta+1}\preceq M^*$ be $\aleph_\eps$-prime over 
$M_\beta\cup N_{\nu}$.

Note that in either case, $R^\infty(c^*/M_{\beta+1})<R^\infty(c^*/M_{\beta})$, so by continuing in this fashion,
$c^*$ will be contained in $M_{\beta+k}$ for some finite $k$.  
\endproof

Suppose that $M\preceq M^*$, and that $\d^*$ is a prime $\aeP$-decomposition of $M^*$ respecting $M$, as witnessed
by the sequences $\<M_\alpha\>$, $\<\d_\alpha\>$, $\<a_\alpha\>$.
Let $(\star)_\alpha$ denote the statement:
\begin{quotation}
For all finite sets $A\subseteq M$, $B\subseteq M_\alpha$, and finite subtree $t\subseteq I_\alpha$,
there is a finite set $A^*\subseteq M$ containing $A$ and a finite subtree $t^*\subseteq I_\alpha$
containing $t$ such that $\tp(B/\bigcup\{N_\rho:\rho\in t^*\})$ is $\aleph_\eps$-isolated and
$\fg M {A^*} {\{N_\rho:\rho\in t^*\}B}$.
\end{quotation}

\begin{Lemma}   $(\star)_\alpha$ holds for all $\alpha\le\alpha^*$.
\end{Lemma}

\bp  We prove this by induction on $\alpha$.
For $\alpha=0$, this is immediate since $M_0=N_{\<\>}$ and is independent from $M$ over $\emptyset$,
hence over any finite subset of $M$.
For $\alpha$ a non-zero limit ordinal, this follows easily from superstability.

For the successor case, fix $\alpha=\beta+1$ and assume that $(\star)_\beta$ holds.
The verification of $(\star)_\alpha$ splits into two cases, depending on whether or not
$I_\beta$ is extended.  Here, we discuss the case where $I_\alpha=I_\beta\cup\{\nu\}$
and leave the other (easier) case to the reader.  So $N_\nu$ is $\aleph_\eps$-prime
over $N_{\nu^-}\cup\{a_\beta\}$, $\fg {N_\nu} {N_{\nu^-}} {M_\beta}$, and $M_\alpha$
is $\aleph_\eps$-prime over both sets $\bigcup\{N_\rho:\rho\in I_\alpha\}$ and $M_\beta\cup N_\nu$.

Towards verifying $(\star)_\alpha$, fix finite sets $A\subseteq M$, $B\subseteq M_\alpha$,
and a finite subtree $t\subseteq I_\alpha$.  Begin by choosing finite sets $C_\nu\subseteq N_\nu$
and $C_\beta\subseteq M_\beta$ such that 
$$\stp(B/C_\beta C_\nu)\vdash \stp(B/M_\beta N_\nu)$$
Without loss, we may assume $a_\beta\in C_\nu$ and $C_\eta\cup C_\beta\subseteq B$.

Next, by superstability choose finite sets $D\subseteq\bigcup\{N_\rho:\rho\in I_\beta\}$ and
$A'\subseteq M$ containing $A$ such that
$$\fg {C_\nu} {DA'} {\bigcup\{N_\rho:\rho\in I_\beta\} M}$$.
Similarly, choose finite sets $E_\beta\subseteq M_\beta$ and $A''\subseteq M$ containing $A'$
such that $$\fg B {E_\beta A''} {M_\beta M}$$
Without loss, we may assume $D\subseteq E_\beta$, $\nu\in t$, and $D\subseteq \bigcup\{N_\rho:\rho\in s\}$,
where $s:=t\setminus\{\nu\}$.

Now apply $(\star)_\beta$ to the triple $(A'', E_\beta,s)$ and get a finite set $A^*\subseteq M$ and 
a finite tree $s^*\subseteq I_\beta$.  Let $t^*:=s^*\cup\{\nu\}$.  We claim that $(A^*,t^*)$ are as desired
in the statement of $(\star)_\alpha$.

\medskip

\par\noindent{\bf Claim 1:}  $B/\bigcup\{N_\rho:\rho\in t^*\}$ is $\aleph_\eps$-isolated.

\medskip  To see this, first note that $C_\beta\subseteq M_\beta$ is $\aleph_\eps$-isolated over
$\bigcup\{N_\rho:\rho\in s^*\}$.  Since $\fg {M_\beta} {N_{\nu^-}} {N_\nu}$ and $N_{\nu^-}$ is
$\aleph_\eps$-saturated, it follows that $C_\beta$ is $\aleph_\eps$-isolated over $\bigcup\{N_\rho:\rho\in t^*\}$
as well.  Also, $C_\nu\subseteq N_\nu$, so it follows immediately that
$C_\beta C_\nu/\bigcup\{N_\rho:\rho\in t^*\}$ is $\aleph_\eps$-isolated as well.
But, as $\stp(B/C_\beta C_\nu)\vdash \stp(B/\bigcup\{N_\rho:\rho\in t^*\}$, the result follows.

\medskip

\par\noindent{\bf Claim 2:}  $\fg M {A^*} {N_0 N_\nu B}$, where $N_0:=\bigcup\{N_\rho:\rho\in s^*\}$.

\medskip
First, it follows from our application of $(\star)_\beta$ that
$\fg {M} {A^*} {N_0 E_\beta}$.  We next consider $C_\nu$.
By the definition of $E_\beta$ and $A''$ we have $\fg {C_\nu} {E_\beta A''} {M_\beta M}$.
So, by monotonicity, we have $\fg {C_\nu} {E_\beta A^*} {N_0 M}$, hence
$\fg {C_\nu} {E_\beta A^* N_0} {M}$.  Thus,  the transitivity of non-forking yields
$$\fg M {A^*} {N_0 E_\beta C_\nu}$$
Finally,  our choice of $N_\nu$ gives $\fg {N_\nu} {N_{\nu^-} a_\beta} {M_\beta M}$.
But $a_\beta\in C_\nu\subseteq N_\nu$, so
$\fg {N_\nu} {N_{\nu^-} C_\nu} {N_0 E_\beta A^* M}$.  As $N_{\nu^-}\subseteq N_0$,
monotonicity yields
$$\fg M {N_0 E_\beta A^*} {N_\nu}$$
and we finish by quoting the transitivity of non-forking.
\endproof

\begin{Proposition}  \label{finiteerror}
Suppose that  $M\preceq M^*$ with  $M^*$ saturated and $||M^*||>||M||+2^{|T|}$.
If $\d^*$ is a prime $\aeP$-decomposition of $M^*$ respecting $M$, then for every finite
$A\subseteq M$ and every finite subtree $t\subseteq I_{\d^*}$, there is a finite set $A^*\subseteq M$
containing $A$, a finite subtree $t^*\subseteq I_{\d^*}$ extending $t$, 
and $M_{t^*}\preceq M^*$ that is $\aleph_\eps$-prime over $\bigcup\{N_\rho:\rho\in t^*\}$
such that $A\subseteq M_{t^*}$, but $\fg M {A^*} {M_{t^*}}$.
\end{Proposition}

\bp  Fix finite $A\subseteq M$ and $t\subseteq I_{\d^*}$.  If $M^*=M_{\alpha^*}$, then 
applying $(\star)_{\alpha^*}$  to the triple $(A,A,t)$ yields a finite set $A^*\subseteq M$
containing $A$ and $t^*$ such that $\tp(A/\bigcup\{N_\rho:\rho\in t^*\})$ is $\aleph_\eps$-isolated
and $\fg M {A^*} {\{N_\rho:\rho\in t^*\}}$.  Thus, as $M^*$ is saturated, we can find
$M_{t^*}\preceq M^*$ containing $A$ that is both $\aleph_\eps$-prime over $\bigcup\{N_\rho:\rho\in t^*\}$
and is independent from $M$ over $A^*$.
\endproof

\subsection{A weak uniqueness theorem for $\P$-decompositions}

The goal of this subsection is Theorem~\ref{quotable}, which is used in \cite{LSh2}.
As we only seek a sufficient condition, the statements and assumptions in Theorem~\ref{quotable}
are inelegant at best.  Additionally, throughout this subsection we assume

\centerline{{\bf  $T$ is totally transcendental with $\P$-NDOP and $\P=\Pactive$}}
The assumption of the theory $T$ being totally transcendental is only used in
Lemma~\ref{freemodel}, and one could easily imagine it being replaced by much
weaker assumptions.
We begin with a standard fact about superstable theories.

\begin{Lemma}  \label{stable}
Suppose that $p\in S(A)$ is stationary and that $J$
is an infinite, $A$-independent set of realizations of $p$.  Let $B\supseteq
A\cup J$, let $p'\in S(B)$ denote the non-forking extension of $p$,
and let $C\supseteq B$ be constructible over $B$.
Then $p'$ has a unique extension to $S(C)$.
\end{Lemma}

\begin{Definition}  {\em Given any model $M$, a {\em $\Pr$-decomposition 
$\d=\<M_\eta,a_\eta:\eta\in I\>$ inside $M$\/} is a weak $\P$-decomposition
inside $M$ with the additional property that $\tp(a_\nu/M_{\nu^-})\in\P$ (hence is
regular) for every $\nu\in I\setminus\{\<\>\}$.
$\d$ is a {\em $\Pr$-decomposition of $M$\/} if, in addition, for every $\eta\in I$,
$\{a_\nu:\nu\in Succ(\eta)\}$ is a maximal $M_\eta$-independent set of realizations
of types in $\P$.
A $\Pr$-decomposition of $M$ is {\em $\P$-finitely saturated\/} if, for every
$\epsilon$-finite $A\subseteq M$ and $b\in M$ such that 
$\tp(b/A)\in\P$, there is some $\eta\in I$
such that $\tp(b/A)\not\perp M_\eta$.
}
\end{Definition}

As notation, given a $\Pr$-decomposition $\d=\<M_\eta,a_\eta:\eta\in I\>$ of $M$,
let $I'=I\setminus\{\<\>\}$.  For each $\eta\in I'$, let
$p_\eta=\tp(a_\eta/M_{\eta^-})$ and fix an $\epsilon$-finite 
$A_\eta\subseteq M_{\eta^-}$ over which $p_\eta$ is based and stationary.
We let $\PP{{a_\eta}\choose A_{\eta}}$ abbreviate $\PP({{\acl(A_\eta a_\eta)}\choose {A_\eta}},\C)$.
Note that by Proposition~\ref{compare}(1),
$\PP{{a_\eta}\choose A_{\eta}}=
\PP{{a_\eta}\choose A_{\eta}'}$ for any $\epsilon$-finite $A'_\eta\subseteq M_{\eta^-}$ on which $p_\eta$ is
based and stationary.

Let $C_\eta:=\{\rho\in I':\rho^-=\eta^-$ and $p_\rho=p_\eta\}$
and let $J_\eta:=\{a_\rho:\rho\in C_\eta\}$.

\begin{Lemma}  \label{freemodel}
Fix any $\Pr$-decomposition $\d=\<M_\eta,a_\eta:\eta\in I\>$
of $M$ and choose any $\eta\in I'$ for which $C_\eta$ is infinite.  Denote
$p_\eta,A_\eta,C_\eta,J_\eta$ by $p,A,C,J$, respectively.  
For any $b\in\C$ realizing $p|_A$, if $\fg b A {M_{\eta^-}J}$ then
$\fg b A M$.
\end{Lemma}

\bp  Fix any element $b$ such that
$\fg b A {M_{\eta^-}}J$.  Let $D:=\bigcup\{M_\rho:\rho\in C\}$ and let 
$E:=\bigcup\{M_\nu:\nu\in I\}$.  First, as $J$ is infinite,
$$\tp(b/M_{\eta^-}J)\vdash\tp(b/D)$$ by Lemma~\ref{stable}.
Next, $\tp(b/D)\vdash\tp(b/E)$ by the independence of the tree, orthogonality,
and the non-forking calculus.  
Next, form a maximal, continuous elementary chain of submodels $\<M_\alpha:\alpha<\beta\>$
of $M$ such that $M_0$ is constructible over $E$, and given $M_\alpha$, $M_{\alpha+1}$ is
constructible over $M_\alpha\cup\{b_\alpha\}$ for some $b_\alpha$ such that $\tp(b_\alpha/M_\alpha)$
is regular.  (Here is where we use the assumption that $T$ is totally transcendental.)  
Clearly, the maximality of the sequence implies that the union is all of $M$.
However, by Lemma~\ref{stable} and the fact that $\tp(b_\alpha/M_\alpha)\perp\P$ (which follows
from $\P=\Pactive$) we conclude that
$$\tp(b/E)\vdash\tp(b/M)$$ 
That $\fg b A M$ follows by the transitivity
of non-forking.
\endproof

\begin{Lemma} \label{comparable} Suppose that 
$\d=\<M_\eta,a_\eta:\eta\in I\>$ is a $\Pr$-decomposition of $M$ and there
is $q\in\P$ and $\eta\in\max(I')$ such that $q\not\perp M_\eta$, but $q\perp M_{\eta^-}$.
Then, for any $\nu\in I$,
$$\nu\triangleleft\eta\qquad\hbox{if and only if}\qquad q\in\PP{{a_\nu}\choose{A_\nu}}$$
\end{Lemma}

\bp
First, assume that $\nu\triangleleft\eta$.  
Let $\d_0:=\<M_\delta,a_\delta:\nu^-\trianglelefteq\delta \trianglelefteq\eta\>$.  
As in the proof of Lemma~\ref{blowup}, we can blow up $\d_0$ to
a sequence 
$\d_0^*:=\<M^*_\delta,a_\delta:\nu^-\trianglelefteq\delta \trianglelefteq\eta\>$,
where $\d_0^*$ is an $\aeP$-decomposition inside $\C$, with $q\not\perp M^*_\eta$,
but $q\perp M^*_{\eta^-}$.  Thus, $q\in\PP{{a_\nu}\choose{A_\nu}}$ by its definition
and Lemma~\ref{straightBA}(3).

Conversely, assume by way of contradiction that
$q\in\PP{{a_\nu}\choose{A_\nu}}$ but $\neg(\nu\triangleleft\eta)$.  As $\nu\neq\eta$
and $\eta\in\max(I')$, $\nu$ and $\eta$ are incomparable.  However, since 
$q\in\PP{{a_{\eta^-}}\choose{A_{\eta^-}}}$ from above, it follows from 
Corollary~\ref{incomparable} that $\nu$ and $\eta^-$ are comparable.  Thus,
$\eta^-\triangleleft\nu$.  But then, as $q\perp M_{\eta^-}$ and $\fg  {M_\eta} {M_{\eta^-}} {a_\nu M_{\nu^-}}$,
it follows that $q$ is orthogonal to any chain starting with $M_{\nu^-}$ and
$a_\nu$.   

\begin{Definition}  {\em  Suppose $S\subseteq\P$.  A $\Pr$-decomposition
$\d=\<M_\eta,a_\eta:\eta\in I\>$ (inside $\C$) {\em supports $S$\/}
if, for every $q\in S$, there is a (unique) $\eta(q)\in \max(I')$ such that
$q\not\perp M_{\eta(q)}$, but $q\perp M_{\eta(q)^-}$.
If $\d$ supports $S$, we let
\begin{itemize}
\item  Field$(S):=\{\eta(q)\in\max(I'):q\in S\}$; and
\item $I^S:=\{\nu\in I:\nu\triangleleft\eta$ for some $\eta\in{\rm Field}(S)\}$.
\end{itemize}
}
\end{Definition}

\begin{Lemma}  \label{basicfacts}
Suppose $S\subseteq\P$ and fix a $\Pr$-decomposition
$\d=\<M_\eta,a_\eta:\eta\in I\>$ (inside $\C$) that supports $S$.
Then:
\begin{enumerate}
\item  If $\nu\in I^S$, then $\tp(a_\nu/M_{\nu^-})$ is trivial;
\item  for $\nu\in I'$, $\nu\in I^S$ if and only if
$\PP{{a_\nu}\choose{A_\nu}}\cap S\neq\emptyset$; and
\item  if, for all $\delta\in I^S$, there is a single
$\epsilon$-finite $A^*\subseteq M_\delta$
such that $\tp(a_\nu/M_\delta)$ is based and stationary on $A^*$
for every $\nu\in Succ_{I^S}(\delta)$,
then for any $\nu\in Succ_{I^S}(\delta)$ and any $b\in\C$
realizing $\tp(a_\nu/A^*)$, if $\PP{b\choose{A^*}}\cap S\neq\emptyset$,
then $\fg b {A^*} {M_\delta}$.
\end{enumerate}
\end{Lemma}

\bp  (1)  It follows immediately from the definition of  $\Pr$-decompositions
and $I^S$ that
$\tp(a_\nu/M_{\nu^-})\in \P$ and has positive $\P$-depth.  Hence,
the type is trivial by Lemma~\ref{positivedepth}.

(2)  This is immediate from unpacking the definitions and Lemma~\ref{basicfacts}.

(3)  Choose $A^*,\delta,\nu,$ and $b$ as required.  Choose 
$r\in \PP{b\choose{A^*}}\cap S$ and look at $\eta(r)\in\max(I')$.  
By Lemma~\ref{basicfacts}, $\delta\triangleleft\eta(r)$.
Choose $\mu\in Succ_{I^S}(\delta)$ satisfying $\mu\triangleleft\eta(r)$.
By our choice of $A^*$ and Lemma~\ref{basicfacts} again, $r\in\PP{{a_\mu}\choose{A^*}}$, so by Proposition~\ref{compare}(5),
$\nfg b {A^*} {a_\mu}$  But then, as $\tp(b/A^*)$ is a trivial regular type,
$b$ is domination equivalent to $a_\mu$ over $A^*$.  Since
$\fg {a_\mu} {A^*} {M_\delta}$, we conclude that the same holds for
$b$.
\endproof

\begin{Definition}  {\em Fix 
$S\subseteq \P$ and a model $M$.
A $\Pr$-decomposition $\d=\<M_\eta,a_\eta:\eta\in I\>$ of $M$  is
{\em $S$-reasonable} if 
\begin{enumerate}
\item  $\d$ is $\P$-finitely saturated and supports $S$;
\item for each $\eta\in I'$:
\begin{enumerate}
\item    $C_\eta\cap I_S$ is infinite;
\item  $p_\rho=p_\eta$ iff $p_\rho\not\perp p_\eta$ 
for every $\rho\in I'$ such that $\rho^-=\eta^-$; and
\item  If $b\in\C$ and $\tp(b/A_\eta)=p_\eta|_{A_\eta}$ and 
$\PP{b\choose {A_\eta}}\cap S\neq\emptyset$, then 
$\fg b {A_\eta} {M_{\eta^-}}$.
\end{enumerate}
\end{enumerate}
}
\end{Definition}

\begin{Definition}  {\em  A {\em weak bijection} between two infinite
sets $I$ and $J$ is a bijection $h:I'\rightarrow J'$, where $I',J'$ are
cofinite subsets of $I,J$, respectively.
}
\end{Definition}

As notation, for $\eta\in I^S\setminus\{\<\>\}$, let 
$J_\eta^S=\{a_\rho:\rho\in C_\eta\cap I^S\}$.

\begin{Proposition}  \label{corresponding}
Fix a set $S\subseteq \P$ and a model $M$.
For $\ell=1,2$, let
 $\d_\ell=\<M_{\eta_\ell},a_{\eta_\ell}:\eta_\ell\in I_\ell\>$
be two $S$-reasonable $\Pr$-decompositions of $M$. 
For any $\eta_\ell\in I^S_\ell$, choose $\eta_{3-\ell}\in I_{3-\ell}$
such that
$p_{\eta_1}\not\perp p_{\eta_2}$.
There is a weak bijection $h:J_{\eta_1}^S\rightarrow J_{\eta_2}^S$
satisfying $\PP{a\choose{A_{\eta_1}}}=\PP{{h(a)}\choose{A_{\eta_2}}}$
for each $a\in\dom(J_{\eta_1})$.
\end{Proposition}

\bp  For definiteness,
assume we have that $\eta_1\in I^S_1$.
 Let $E=A_{\eta_1}\cup A_{\eta_2}$.  For $\ell=1,2$,
let $p_\ell\in S(E)$ be parallel to $p_{\eta_\ell}$, let
$J'_\ell=\{a\in J_{\eta_\ell}:\fg a {A_{\eta_\ell}} E\}$,  and  let
$$J^S_\ell=\{a\in J'_\ell:\PP{a\choose{A_{\eta_\ell}}}\cap S\neq\emptyset\},$$
which is a cofinite subset of $J^S_{\eta_\ell}$.  In particular,
$J^S_1\neq\emptyset$ since $C_{\eta_1}\cap I^S$ is infinite.
As well, choose
 a maximal $E$-independent set  $J^*_\ell$
of realizations
of $p_{\ell}$ in $\C$ extending $J_\ell$.
As  $p_1$ and $p_2$ are non-orthogonal
trivial regular types, it follows from Proposition~\ref{trivialapp}
that  there is a unique
bijection $h:J_1^*\rightarrow J_2^*$ satisfying $\nfg {h(a)} E a$ for each
$a\in J_1^*$.

As $\fg {a_\ell} {A_{\eta_\ell}} E$ for $\ell=1,2$ and every
$a_\ell\in J^*_\ell$,  by Proposition~\ref{compare}(1) we have that
$$\PP{a\choose{A_{\eta_1}}}=\PP{a\choose E}=\PP{{h(a)}\choose E}=
\PP{{h(a)}\choose{A_{\eta_2}}}$$
for each $a\in J^*_\ell$.

\medskip\par\noindent{\bf Claim.} For every $a\in J^S_1$, $h(a)\in J^S_2$.

\medskip\noindent
{\bf Proof of Claim:}  Choose any $a\in J^S_1$.  
We first find an element $b\in J_{\eta_2}$ such that
$\nfg {h(a)} E b$.  Since $a=a_\rho$ for some
$\rho\in I^S_1$ satisfying $\rho^-=\eta_1^-$, 
$\PP{a\choose{A_{\eta_1}}}\cap S\neq\emptyset$.  As the two sets are equal,
$\PP{{h(a)}\choose{A_{\eta_2}}}\cap S\neq\emptyset$ as well.
As $\d_2$ is $S$-reasonable, this implies $\fg {h(a)} {A_{\eta_2}} {M_{\eta_2^-}}$
Next, we argue that $h(a)$ must fork with $J_{\eta_2}$ over $M_{\eta_2^-}$,
because if this were not the case, then by Lemma~\ref{freemodel} we would
have $\fg {h(a)} E M$.  But, as $\nfg a E {h(a)}$, the fact that $p_{\eta_2}$
has weight one would imply that $\fg a E M$, which is absurd since $a\in M$.

Thus, $h(a)$  forks with $J_{\eta_2}$ over $M_{\eta_2^-}$.
By triviality, there is a unique $b\in J_{\eta_2}$ such that $\nfg {h(a)} {M_{\eta_2^-}}
b$.  However, as both $h(a)$ and $b$ are free from $M_{\eta_2^-}$
over $A_{\eta_2}$, it follows that $h(a)$ and $b$ fork over $A_{\eta_2}$,
completing the first part of our argument.  
 
Next, since $h(a)$ realizes $p_2$, it is free from $E$ over $A_{\eta_2}$.  As $p_{\eta_2}$
has weight one, the last two statements imply that $b$ is free from $E$
over $A_{\eta_2}$ as well.  Thus, $b\in J_2'$.  As well, we have
that $\PP{a\choose E}=\PP{b\choose E}$, so the latter has non-empty intersection
with $S$.  Thus, $b\in J^S_2$.

Finally, note that both $h(a)$ and $b$ are elements of $J^*_2$ that fork
with each other over $E$.  Thus, $h(a)=b$ by the $E$-independence of   $J_2^*$.  So $h(a)\in J^S_2$, completing the proof of the Claim.  
\endproof

It follows from the Claim that $J_2^S$ is non-empty.  Once we know this,
the situation becomes symmetric, so by running the Claim backwards,
$h^{-1}$ maps $J^S_2$ into $J^S_1$.  That is, the restriction of
$h$ to $J^S_1$ is a bijection with $J^S_2$, which completes the proof of
the Proposition.
\endproof

We set some notation about partial maps between trees.
Given a tree $I$, a {\em large subtree} of $I$ is a non-empty (downward closed)
subtree $J$ such that for every $\eta\in J$, $Succ_I(\eta)\setminus J$
is finite.  Given two trees $J$ and $K$, an {\em almost embedding $h$ from $J$ to $K$}
has $\dom(h)$ a large subtree of $J$, $range(h)\subseteq K$, 
$h(\<\>_J)=\<\>_K$, and for all $\eta,\nu\in\dom(h)$,
$$\eta\triangleleft \nu \quad\hbox{if and only if}\quad h(\eta)\triangleleft h(\nu)$$
The trees $J$ and $K$ are {\em almost isomorphic\/} if there is an almost embedding
$h$ from $J$ to $K$ in which $range(h)$ is a large subtree of $K$.

For $J$ any tree and $\nu\in J$, let $J_{\trianglerighteq\nu}$ be the tree with root
$\nu$ and universe $\{\eta\in J:\eta\trianglerighteq\nu\}$.
Given two trees $J$ and $K$ and $\nu\in J$, $\mu\in K$, an 
{\em almost embedding $h$ from $J$ to $K$ over ($\nu,\mu)$}
is an almost embedding from $J_{\trianglerighteq\nu}$ to $K_{\trianglerighteq\mu}$.

Finally, if $J$ and $K$ are trees indexing  decompositions, we call a pair
$(\eta,\nu)\in J\times K$ {\em $\PP$-equivalent} if
either $\eta=\<\>=\nu$, or both $\eta,\nu\neq\<\>$ and 
$\PP{{a_\eta}\choose{A_\eta}}=\PP{{a_{\nu}}\choose{A_{\nu}}}$.
An {\em almost $\PP$-embedding from $J$ to $K$} is an almost embedding $h$
from $J$ to $K$ with the pair $(\eta,h(\eta))$ $\PP$-equivalent for each
$\eta\in\dom(h)$.
Note that if $h$ is an almost $\PP$-embedding and $h(\eta)=\nu$, then the
restriction of $h$ to $J_{\trianglerighteq\eta}:=\{\delta\in\dom(h):\delta
\trianglerighteq\eta\}$ is an almost $\PP$-embedding over $(\eta,\nu)$.

Given all of this notation, the proof of the following Corollary simply involves successively
iterating Proposition~\ref{corresponding}, using the fact that each decomposition
is $\P$-finitely saturated.

\begin{Corollary} \label{almostembedding} 
Fix a set $S\subseteq \P$ and a model $M$.
For $\ell=1,2$, suppose that $\d_\ell=\<M_{\eta_\ell},a_{\eta_{\ell}}:\eta_\ell\in I_\ell\>$
are $S$-reasonable $\Pr$-decompositions of $M$ with the additional property that for each $\ell$ and $\nu_\ell\in I_\ell$,
$$\{p:\hbox{there is}\ \eta_\ell\in Succ(\nu_\ell)\ \hbox{such that}\ p_{\eta_\ell}=p\wedge\PP{{a_{\eta_\ell}}
\choose{A_{\eta_\ell}}}\cap S\neq\emptyset\}$$
is finite.  Then:
\begin{enumerate}
\item  For $\ell=1,2$, there is an almost $\PP$-embedding $h$ from $I^S_\ell$ to $I^S_{3-\ell}$; and
\item  For $\ell=1,2$ and any $\PP$-equivalent pair 
$(\eta_\ell,\eta_{3-\ell})\in I^S_\ell\times I^S_{3-\ell}$ there is
an almost $\P$-embedding from $I^S_\ell$ to $I^S_{3-\ell}$ over $(\eta_\ell,\eta_{3-\ell})$.
\end{enumerate}
\end{Corollary}

If we wish to conclude more, namely that the trees $I^S_1$ and $I^S_2$ are almost isomorphic, then we
need show that the almost embeddings given above preserve lengths, i.e.,
that $\lg(h(\eta))=\lg(\eta)$ for every $\eta\in\dom(h)$.  
To accomplish this,
we need to put additional constraints on the shapes of the trees $I^S$.
The conditions we require are severe, but will be easily satisfied in our
construction in \cite{LSh2}. 

\begin{Definition} \label{proper} {\em  
A {\em two-coloring\/} of a tree $I$ is a sequence $\<E_\eta:\eta\in I\>$
where each $E_\eta$ is an equivalence relation on $Succ(\eta)$ with
at most two classes, each of which is infinite.  (If $Succ(\eta)=\emptyset$,
then of course $E_\eta$ is empty as well.)
 A node $\eta\in I$ has {\em uniform depth $n$\/} if every branch of the tree 
$I_{\trianglerighteq \eta}$ has length exactly $n$.  A node $\eta$ {\em often has unbounded depth\/}
if every large subtree $J\subseteq I_{\trianglerighteq\eta}$ has an infinite branch.
A node $\eta$ is an {\em $(m,n)$-cusp} if there are infinite sets $A_m,A_n,B\subseteq Succ(\eta)$ such that
\begin{enumerate}
\item the set $A_m\cup A_n$ is pairwise $E_\eta$-equivalent;
\item  each $\delta\in A_m$ has uniform depth $m$;  
\item  each $\rho\in A_n$ has uniform depth $n$;
and
\item  each $\gamma\in B$ is often unbounded.
\end{enumerate}
A {\em cusp\/} is an $(m,n)$-cusp for some $m\neq n$.

Fix any function $\Phi:\omega\rightarrow\omega$.
We say the two-colored
tree $I$ is {\em $\Phi$-proper\/} if, for every node $\eta\in I$, 
\begin{enumerate}

\item  either $\eta$ has uniform depth $n$ for
some $n$, or else $\eta$ often has unbounded depth; 
\item  if $\eta$ is an $(m,n)$-cusp, then $\lg(\eta)=\Phi(m-n)$;
\item  if $E_\eta$ has two classes, then $\eta$ is a cusp;
\item  if $J$ is a large subtree of $I$, $\eta\in J$ is often unbounded, then there is a cusp $\nu\in J$ with $\nu\trianglerighteq
\eta$.

\end{enumerate}
}\end{Definition}

Note that if $I$ is a two-colored tree satisfying the conditions above,
then for every $\gamma\in I$ that is of any uniform depth $k$,
there are a unique $\eta,\delta$ satisfying $\delta\trianglelefteq\gamma$,
$\eta=\delta^-$, $\eta$ is a cusp, and $\delta$ has uniform depth $n$ for
some $n\ge k$.

\begin{Lemma} \label{match}
 Suppose that $M,S,\d_1,\d_2$ satisfy the assumptions
of Corollary~\ref{almostembedding} and additionally assume that
both $I^S_1, I^S_2$, when two-colored by the relations
$E_\eta$ defined by $E_\eta(\delta,\rho)$ iff $\delta^-=\eta=\rho^-$
and $p_\delta=p_\rho$, are $\Phi$-proper for the same  function $\Phi$.
Then for every $\PP$-equivalent pair 
$(\eta,\nu)\in I^S_1\times I^S_2$,
\begin{enumerate}
\item  $\eta$ is often unbounded in $I^S_1$ if and only if $\nu$
is often unbounded in $I^S_2$;
\item  for any $n$, $\eta$ has uniform depth $n$ if and only if
$\nu$ has uniform depth $n$;
\item  if $\lg(\eta)=\lg(\nu)$ and $\eta$ has uniform depth $n$ for some $n$,
then any almost $\PP$-embedding over $(\eta,\nu)$ preserves lengths; and
\item  if $\lg(\eta)\le\lg(\nu)$ and $\eta$ is an $(m,n)$-cusp, then
$\nu$ is also an $(m,n)$-cusp, $\lg(\eta)=\lg(\nu)$, and for any almost
$\PP$-embedding $h$ over $(\eta,\nu)$, $\lg(h(\delta))=\lg(\delta)$ for
all $\delta\in\dom(h)\cap Succ(\eta)$ of uniform depth $m$ or $n$;
\item  if $\lg(\eta)=\lg(\nu)$ then every almost $\PP$-embedding over 
$(\eta,\nu)$ preserves lengths; and
\item if $\lg(\eta)=\lg(\nu)$, then the number of $E_\eta$-classes in $I^S_1$
equals the number of $E_\nu$-classes in $I^S_2$.
\end{enumerate}
\end{Lemma}

\bp  (1) First assume that $\eta$ is often
unbounded.  By Corollary~\ref{almostembedding}(2), choose 
an almost $\PP$-embedding $h$ from $I^S_1$ to $I^S_2$ over $(\eta,\nu)$.
Choose a strictly $\triangleleft$-increasing sequence $\<\eta_n:n\in\omega\>$
from $\dom(h)$ with $\eta_0=\eta$.  Then $\<h(\eta_n):n\in\omega\>$
is a strictly $\triangleleft$-increasing sequence in $I^S_2$ with
$h(\eta_0)=\nu$.  Thus, $\nu$ cannot have any finite uniform depth,
so it must be often unbounded by properness.  The converse is symmetric.

(2)  Suppose that $\nu$ has uniform depth $n$.  Then by (1), $\eta$ has uniform
depth $m$ for some $m$.  Arguing as in (1), $m\le n$, since 
if we choose any almost $\PP$-embedding $h$ from $I^S_1$ to $I^S_2$ 
over $(\eta,\nu)$, then the image of any strictly $\triangleleft$-increasing
sequence $\<\eta_i:i<m\>$ with $\eta_0=\eta$ would be a strictly
$\triangleleft$-increasing sequence of length $m$ over $\nu$.
But then, by symmetry, we would also have $n\le m$, so $n=m$.
The converse is symmetric.

(3)  Suppose that $h$ is any almost $\PP$-embedding
over $(\eta,\nu)$, where $\lg(\eta)=\lg(\nu)$, $\eta$ has uniform depth
$n$.   Then $\nu$ also has uniform depth $n$.  So, every maximal 
$\triangleleft$-increasing sequence extending $\eta$ has length $n$,
the image of any such sequence under $h$ is also a strictly $\triangleleft$-increasing
sequence of length $n$, but there is no strictly $\triangleleft$-increasing
sequence of length more than $n$ extending $\nu$.  Thus, $h$ must map
immediate successors to immediate successors, and consequently preserve
lengths.

(4) Suppose that  $\eta$ is an $(m,n)$-cusp and $\lg(\eta)\le\lg(\nu)$.  
Choose an almost 
$\PP$-embedding $h$ from $I^S_1$ to $I^S_2$ over $(\eta,\nu)$.
Choose $E_\eta$-equivalent
$\delta\in Succ(\eta)\cap\dom (h)$ of uniform depth $m$
and $\rho\in Succ(\eta)\cap\dom(h)$ of uniform depth $n$.
Choose $\mu\in I^S_2$ and $q\in S(M^2_\mu)$ such that
$p_\delta$ (which $=p_\rho$) is non-orthogonal to $q$.
By the definition of $h$, both $h(\delta),h(\rho)\in Succ(\mu)$.
We argue that $\mu=h(\eta)$.  To see this, first note
that since $h$ is $\triangleleft$-preserving, $h(\eta)\triangleleft h(\delta)$
and $h(\eta)\triangleleft h(\rho)$, so $h(\eta)\trianglelefteq \mu$.
But, it follows from (2) 
that $h(\delta)$ is uniformly of depth $m$ and $h(\rho)$ is uniformly
of depth $n$.  Thus, $\mu$ is an $(m,n)$-cusp and hence
$\lg(\mu)=\Phi(m-n)=\lg(\eta)$.  As we assumed that
$\lg(\eta)\le\lg(\nu)$ and $h(\eta)=\nu$, we have that $\lg(\mu)=\lg(h(\eta))$,
hence $\mu=h(\eta)=\nu$.  This yields $\lg(\nu)=\lg(\eta)$.
Finally, the argument above showed that $h(\delta)\in Succ(\nu)$ whenever
$\delta\in\dom(h)\cap Succ(\eta)$ has uniform depth $m$ or $n$.

(5)  Assume that $\lg(\eta)=\lg(\nu)$
and fix any almost $\PP$-embedding $h$ from $I^S_1$ to $I^S_2$
over $(\eta,\nu)$.    
Note that $\lg(h(\mu))\ge\lg(\mu)$ for any $\mu\in\dom(h)$
simply because $h$ is $\triangleleft$-preserving.
We first consider the often unbounded nodes $\mu\in\dom(h)$.
Specifically, we argue by induction on $k$ that $\lg(h(\mu))=\lg(\mu)$
for every often unbounded node $\mu\in\dom(h)$ for which there
is a cusp $\zeta\trianglerighteq\mu$ with $\zeta\in\dom(h)$
and $\lg(\zeta)=\lg(\mu)+k$.  

When $k=0$, this means that any such $\mu$ is itself a cusp, so
$\lg(h(\mu))=\lg(\mu)$ by (4).  
Next, assume that the statement holds for $k$, and choose $\mu\in\dom(h)$
with some cusp $\zeta\in\dom(h)$ with $\mu\trianglelefteq\zeta$ and
$\lg(\zeta)=\lg(\mu)+k+1$.  Choose $\rho\in Succ(\mu)$ with
$\mu\trianglelefteq\rho\trianglelefteq\zeta$.  Then $\lg(h(\rho))=\lg(\rho)$
by our inductive assumption, so $h(\rho)\in Succ(h(\mu))$,
hence $\lg(h(\mu))=\lg(\mu)$ as well.
Thus, we have shown that lengths are preserved for all
often unbounded nodes $\mu\in\dom(h)$.

Next, assume that $\gamma\in\dom(h)$ has uniform depth.  By the remark
following Definition~\ref{proper}, choose $\mu$ and $\delta$ such
that $\mu$ is a cusp, $\mu=\delta^-$, $\delta\trianglelefteq\gamma$, 
and $\delta$ has uniform depth $n$ for some $n\ge k$.
The last sentence of (4) implies that  $\lg(h(\delta))=\lg(\delta)$.
Thus, $\lg(h(\gamma))=\lg(\gamma)$ follows from (3).
So  $h$ is length-preserving.  

(6)  As the hypotheses are symmetric, it suffices to prove that
the number of $E_\eta$-classes is at most the number of $E_\nu$-classes.
Using Corollary~\ref{almostembedding}, choose an
 almost $\PP$-embedding $h$ over $(\eta,\nu)$.
By (5), $h$ maps immediate successors of $\eta$ to immediate
successors of $\nu$.  As well, for each $\delta\in \dom(h)\cap Succ(\eta)$,
$p_\delta\not\perp p_{h(\delta)}$.  As non-orthogonality is an
equivalence relation on regular types, this implies that $h$ maps
$E_\eta$-classes to $E_\nu$-classes, and maps distinct $E_\eta$-classes
to distinct $E_\nu$-classes.  
As  there are at most two $E_\eta$-classes, the inequality follows.
\endproof

\begin{Theorem}  \label{prettygood}
Fix a set $S\subseteq \P$ and a model $M$.
For $\ell=1,2$, suppose that $\d_\ell=\<M_{\eta_\ell},a_{\eta_{\ell}}:\eta\in I_\ell\>$
satisfy the hypotheses of Lemma~\ref{match}.
  Then there is an almost $\PP$-isomorphism $h$ from $I^S_1$ to $I^S_2$.
\end{Theorem}

\bp  Using Corollary~\ref{almostembedding}, choose any almost
$\PP$-embedding $h$ of $I^S_1$ to $I^S_2$ such that, for any
$\delta\in \dom(h)$, $\dom(h)\cap C_\delta$ is a cofinite subset of
$C_\delta$ and ${\rm range}(h)\cap C_{h(\delta)}$ is a cofinite subset
of $C_{h(\delta)}$.  
From Lemma~\ref{match} we know that $h$ preserves levels and,
for each node $\eta\in\dom(h)$, the number of $E_{h(\eta)}$-classes
is equal to the number of $E_\eta$-classes.  It follows that ${\rm range}(h)$
is a large subtree of $I^S_2$, so $h$ is an almost $\PP$-isomorphism
between $I^S_1$ and $I^S_2$.
\endproof

Finally, we exhibit an extreme case, whose hypotheses are satisfied in
\cite{LSh2}.

\begin{Definition}  {\em
Fix $S\subseteq \P$, a model $M$, and a function $\Phi:\omega\rightarrow\omega$.
A $\Pr$-decomposition $\d=\<M_\eta,a_\eta:\eta\in I\>$ of $M$ is {\em
$(S,\Phi)$-simple\/} if
\begin{enumerate}
\item  $\d$ supports $S$ and $\P$-finitely saturates $M$;
\item  for every $\eta\in I^S$
\begin{enumerate}
\item  $Succ_{I^S}(\eta)$ is empty or infinite, but $E_\eta$ is trivial,
i.e., $p_\nu=p_\mu$ for all $\nu,\mu\in Succ_{I^S}(\eta)$;
\item  $\eta$ is either of some finite uniform depth or is a cusp; and
\item  if $\eta$ is an $(m,n)$-cusp, then $\Phi(m-n)=\lg(\eta)$.
\end{enumerate}
\end{enumerate}
}
\end{Definition}

\begin{Theorem}  \label{quotable}
Fix a set $S\subseteq\P$ and a model $M$, and a function 
$\Phi:\omega\rightarrow\omega$.  If $\d_1$ and $\d_2$
are both $(S,\Phi)$-simple $\Pr$-decompositions of $M$,
then the trees $I^S_1$ and $I^S_2$ are almost $\PP$-isomorphic.
\end{Theorem}

\bp  Because of Theorem~\ref{prettygood}, we only need to
verify that the hypotheses of Lemma~\ref{match} are satisfied for
each of the decompositions.
But this is routine, once one notes that Clause~2(b) is satisfied because 
of the triviality of $E_\eta$ and Lemma~\ref{basicfacts}(3).
\endproof

\end{document}